\documentclass[12pt]{article}
\usepackage[top=0.5in, right=0.4in, left=0.4in, bottom=0.8in]{geometry}
\usepackage{amsmath, amsfonts, amsthm, amssymb, enumerate, color, authblk, color, bm, graphicx, multirow, comment, url, xspace, natbib, titlesec, hyperref, ulem, algpseudocode}
\usepackage{amsmath, amsfonts, amsthm, amssymb, enumerate, color, authblk, color, bm, graphicx, multirow, comment, url, xspace, natbib, titlesec, hyperref, ulem, lscape,float}
\usepackage{hhline}
\usepackage{graphicx}
\usepackage{makecell}
\usepackage{cleveref}
\usepackage{lscape}
\usepackage{pdflscape}
\usepackage{longtable}
\usepackage{rotating}
\usepackage{graphicx}
\usepackage{adjustbox}
\usepackage{ upgreek } 
\usepackage{fullpage}
\usepackage{times}
\usepackage{appendix}
\usepackage[linesnumbered,ruled,vlined]{algorithm2e}
\usepackage{hyperref}

\newtheorem{theorem}{Theorem}[]
\newtheorem{lemma}[theorem]{Lemma}
\newtheorem{proposition}[theorem]{Proposition}

\newcommand{\ie}{{\it i.e., }}
\newcommand{\eg}{{\it e.g., }}

\newcommand{\comments}[1]{}
\newcommand{\removed}[1]{}

\title{A Double-oracle, Logic-based Benders decomposition approach to solve the $K$-adaptability problem}
\author[1]{Alireza Ghahtarani\thanks{Dalhousie University, Halifax, Nova Scotia, B3J 1B6, Canada (alireza.ghahtarani@dal.ca).}}
\affil[1]{Department of Industrial Engineering, Dalhousie University}
\affil[2]{Department of Decision Sciences, HEC Montréal}
\author[1]{Ahmed Saif}
\author[1]{Alireza Ghasemi}
\author[2]{Erick Delage \thanks{Department of Decision Sciences, HEC Montréal, Montréal, Québec, H3T 2A7, Canada.}}
\date{}
\begin{document}
\maketitle

\abstract{We propose a novel approach to solve $K$-adaptability problems with convex objective and constraints and integer first-stage decisions. A logic-based Benders decomposition is applied to handle the first-stage decisions in a master problem, thus the sub-problem becomes a min-max-min robust combinatorial optimization problem that is solved via a double-oracle algorithm that iteratively generates adverse scenarios and recourse decisions and assigns scenarios to $K$-subsets of the decisions by solving $p$-center problems. Extensions of the proposed approach to handle parameter uncertainty in both the first-stage objective and the second-stage constraints are also provided. We show that the proposed algorithm converges to an optimal solution and terminates in finite number of iterations. Numerical results obtained from experiments on benchmark instances of the adaptive shortest path problem, the regular knapsack problem, and a generic $K$-adaptability problem demonstrate the performance advantage of the proposed approach when compared to state-of-the-art methods in the literature.

}

\section{Introduction}\label{sec1}

Robust Optimization (RO) has become a classical framework for dealing with parameter uncertainty in optimization problems \citep{bertsimas2011theory}. In RO, parameter uncertainty is captured through an \emph{uncertainty set} of proper structure and size, and the optimization is conducted with respect to the \emph{worst-case} realization in it. An important class of RO problems that has gained considerable attention recently is adaptive/adjustable robust optimization (ARO), in which some decisions are assumed to be delayable until the realized values of uncertain parameters become partially or fully known. Whereas ARO formulations often lead to better (less pessimistic) solutions than their corresponding static RO models, they are computationally intractable \citep{ben2004adjustable}. However, several exact and approximate algorithms have been proposed to solve important ARO classes such as linear two-stage RO problems \citep{thiele2009robust,chen2009uncertain,kuhn2011primal,zhao2012robust,bertsimas2012adaptive,jiang2012benders,iancu2013supermodularity}. With a few exceptions, these algorithms use duality to handle the second-stage (recourse) problem. Hence, they can be only used with continuous recourse decisions. Although some attempts have been made to develop efficient solution methods for ARO problems with discrete recourse \citep{dhamdhere2005pay,georghiou2015generalized,bertsimas2015design,bertsimas2018binary}, the literature for this class of problems is still sparse.

Recently, an alternative modeling approach, referred to as \emph{$K$-adaptability}, has been proposed as a conservative approximation of ARO problems with discrete recourse \citep{hanasusanto2015k,subramanyam2020k}. Rather than allowing any feasible integer recourse to be selected, the decision-maker \emph{prepares} $K$ solutions in advance (under uncertainty). Then, upon full knowledge of the realized value of the uncertain parameters, the best among these $K$ solutions is selected. Apart from being better in general compared to the solutions of static RO, \cite{buchheim2017min} argued that $K$-adaptability solutions are more easily accepted by a human user as they do not change each time but are taken from a relatively small set of candidate solutions.

Similar to \cite{hanasusanto2015k} and \cite{subramanyam2020k}, we initially focus on a linear version of the problem, in which both the objective function and constraints in the first- and second-stages are affine functions of the decision variables, the uncertain parameters effect the second-stage objective function, and the uncertainty set is polyhedral. However, we later show how the algorithm that we propose can be adapted to solve other variants of the problem. Formally, the linear $K$-adaptability problem under consideration is formulated as follows: 
\begin{subequations}\label{K-adaptability}
\begin{alignat}{3}
&&\min_{\mathrm{x} \in \mathcal{X}, \{\mathrm{y}^{k}\}_{k\in [K]}}&\mathrm{c}'\mathrm{x}+\max_{\mathrm{\upxi} \in \mathrm{\Xi}}\min_{k\in [K]}\left\{\mathrm{\upxi}'\mathrm{Q}\mathrm{y}^{k}:\;\mathrm{T}\mathrm{x}+\mathrm{W}\mathrm{y}^{k}\leq \mathrm{b}\right\}\\
&&\;\text{s.t.}\;\;\;\;\;\; \;& \mathrm{y}^{k}\in \mathcal{Y},\;k\in [K],
\end{alignat}
\end{subequations}
where \(\mathcal{X} \subset \{0,1\}^{n}\) and \(\mathcal{Y} \subset \{0,1\}^{m}\) where their continuous relaxations are polyhedral, \(\Xi \subset \mathbb{R}^{q}\) is a polyhedral set, \(\mathrm{c} \in \mathbb{R}^{n}\), \(\mathrm{Q} \in \mathbb{R}^{q\times m}\), \(\mathrm{T} \in \mathbb{R}^{s\times n}\), \(\mathrm{W} \in \mathbb{R}^{s\times m} \), \(\mathrm{b} \in \mathbb{R}^{s}\), and \([K]=\{1,\dots,K\}\). In this formulation, elements of the vector \(\mathrm{x}\) denote the here-and-now (first-stage) decision variables, \(\xi \in \Xi\) are the uncertain parameters, and \(\mathrm{y}\) are the wait-and-see (second-stage) decision variables. 

A special case of the $K$-adaptability problem, known in the literature as the \emph{min-max-min robust combinatorial optimization} (MMMRCO) problem arises when the only decision that is made under uncertainty is the pre-selection of the $K$ recourse action, \ie the problem does not have the actual first-stage decision. The MMMRCO problem has many practical applications, such as parcel delivery and finding routes in hazardous situations, that are discussed in \cite{arslan2020min}. The basic MMMRCO problem (without constraint uncertainty) is formulated as follows:
\begin{subequations}\label{MMMRCO}
\begin{alignat}{3}
&\;\min_{\mathrm{y}^{k}\in \mathcal{Y}} \max_{\mathrm{\upxi} \in \Xi}\min_{k\in [K]}\{\mathrm{\upxi}'\mathrm{Q}\mathrm{y}^{k}:\;\mathrm{W}\mathrm{y}^{k}\leq \mathrm{b}\}.
\end{alignat}
\end{subequations}

So far, only two solution methods have been developed for the $K$-adaptability problem. \cite{hanasusanto2015k} proposed a mixed-integer linear programming approximation, derived using linear programming
duality, which leads to a monolithic formulation involving bilinear terms. A McCormick envelope is used to linearize the bilinear terms, which requires a large number of new variables to be introduced. Moreover, the number of binary variables increases as $K$ is increased. Consequently, this approach hardly solves the shortest path instances with more than 25 nodes. In another attempt to solve the $K$-adaptability problem, \cite{subramanyam2020k} proposed a branch-and-bound algorithm that enjoys asymptotic convergence in general, but has finite convergence under specific conditions. This algorithm works by generating a relevant subset of uncertainty realizations and enumerating over their assignment to the $K$ solutions. Nevertheless, it is also inefficient for solving shortest path instances with more than 25 nodes.

Besides the aforementioned methods that can handle the general case (\ie with first-stage decisions), a few approaches to solve variants of \eqref{MMMRCO} have been proposed. \cite{chassein2019faster} developed a branch-and-bound algorithm that can solve large instances of the MMMRCO problem, yet only with budget uncertainty sets. They also proposed a heuristic solution algorithm based on the formulation of \cite{hanasusanto2015k}. However, instead of using the McCormick linearization approach to handle the bilinear terms, they used the block-coordinate descent algorithm, which has no optimality guarantee. Moreover, their algorithms hardly solve any instance of the shortest path problem when $K$ grows above 3. \cite{goerigk2020min} developed a heuristic algorithm based on an integer programming formulation and a row-and-column generation algorithm to solve it. Yet, this algorithm has no performance guarantee, which means this algorithm does not guarantee convergence in a finite number of iterations. Recently, \cite{arslan2020min} proposed a solution approach that iteratively generates scenarios of the uncertain parameters and assigns them to solutions by solving a $p$-center problem. However, this algorithm works well only if there is an effective way to restrict and enumerate the search space.

In this paper, we present a new approach to solve the $K$-adaptability problem with binary or integer first-stage decisions. The scenario generation step in the proposed approach enjoys finite convergence when the uncertainty set is polyhedral, but the approach can be used with any convex uncertainty set. Although we focus initially on problems with affine functions, the proposed algorithm can be extended to nonlinear objective function and constraints with respect to the decision variables for second-stage problem. 

The proposed approach uses a \emph{logic-based} (also referred to as \emph{combinatorial}) Benders algorithm to handle the first-stage decisions such that the remaining subproblem is a MMMRCO that is solved iteratively to generate optimality cuts. To solve this subproblem, we propose a double-oracle algorithm that iterates between solving an adversary problem to iteratively generate worst-case scenarios for a $K$-subset of feasible solutions and determine the scenario-solution assignment, and solving a decision-maker's problem to find the optimal $K$-subset of solutions for all the scenarios generated so far. Although the way scenarios are generated and assigned to solutions is similar to that proposed by \cite{arslan2020min}, our approach uses a more efficient way (\ie by solving an optimization problem) to identify the optimal $K$-subset of recourse solutions in every iteration.  

We note that the $K$-adaptability problem formulation provided in \eqref{K-adaptability} is based on that introduced by \cite{bertsimas2010finite}, \ie with a first-stage problem that is not subject to uncertainty, whereas \cite{hanasusanto2015k} addressed an extended version in which both stages are affected by the same uncertain parameters. Hence, we show how the proposed algorithm can be extended to handle $K$-adaptability problems with uncertainty in both stages (whether the two stages depend on the same or different uncertain parameters). Finally, extensive numerical experiments are conducted on benchmark instances of several classical optimization problems, and the computational superiority of the proposed approach \emph{vis-\'{a}-vis} state-of-the-art solution methods is demonstrated.

The remainder of this paper is organized as follows. Section \ref{sec2} presents the approach proposed to solve problem \eqref{K-adaptability} (\ie with linear objective function and constraints and with uncertainty affecting only the recourse objective function). Section \ref{sec3} studies the convergence properties of the proposed algorithm and proves its finite convergence. In section \ref{sec4}, we show how the proposed algorithm can be modified to solve different variants of the $K$-adaptability problem, namely, problems with integer first-stage decision variables, problems with nonlinear functions, and problems affected by uncertainty in the first and second-stages. The numerical experiments conducted to test the proposed algorithm on benchmark problems, and a detailed discussion of their results, are presented in section \ref{sec5}. Finally, section \ref{sec6} provides some conclusions and suggests future research directions.    

\vspace{2mm}
\noindent{\bf Notation.} We use $\mathrm{upright}$ lower and upper case letters, respectively, for vectors and matrices. Individual elements of these vectors and matrices are denoted using $italic$ versions of the same letters. For example, elements of the $J$-dimension vector $\mathrm{x}$ are denoted as $x_j$. Depending on the context, upper case letters might be used also to denote scalars (\eg $J$), whereas lower case letters might denote also functions (\eg $g(\cdot)$). $[J]$ is used as a shorthand for the set of integers $\{1,2,\dots,J\}$ and a partial set of a given set $[J]$ is denoted as $[J']$. We use the symbol $\mathrm{e}$ to denote an all-ones vector of appropriate size. The calligraphic font is used for sets (\eg $\mathcal{X}$). 

\section{The Proposed Solution Approach}\label{sec2}

In this section, we present the proposed approach to solve the $K$-adaptability problem \eqref{K-adaptability} with binary first-stage decision variables and recourse objective uncertainty. First, we show how a logic-based Benders decomposition is applied to deal with the discrete first-stage variables. Then, we describe a double oracles algorithm to solve the subproblem. 

\subsection{A Logic-based Benders Decomposition}

We apply Benders decomposition by projecting the model onto the subspace defined by the first-stage variables $\mathrm{x}$ to get the master problem ($\mathtt{MP}$):
\begin{equation}
\min_{\mathrm{x}\in \mathcal{X}\cap \mathcal{V}}\mathrm{c}'\mathrm{x}+\nu(\mathrm{x}),
\end{equation}
where $\mathcal{V}:=\left\{\mathrm{x}: \mathrm{W}\mathrm{y}^k\leq \mathrm{b}-\mathrm{T}\mathrm{x}~\text{for some}~\mathrm{y}^k\in \mathcal{Y},k\in [K]\right\}$. For a given $\bar{\mathrm{x}}\in \mathcal{X}\cap \mathcal{V}$, $\nu(\bar{\mathrm{x}})$ is the optimal value of the sub-problem ($\mathtt{SP}$)
\begin{subequations}\label{sub-problem}
\begin{alignat}{3}
&\min_{\{\mathrm{y}^{k}\}_{k\in [K]}}\max_{\mathrm{\upxi} \in \Xi}\min_{k\in [K]}\quad&&\mathrm{\upxi}'\mathrm{Q}\mathrm{y}^{k}\\
&\;\;\qquad\qquad\qquad\text{s.t.}\quad &&\mathrm{y}^{k}\in \mathcal{Y},\;\mathrm{W}\mathrm{y}^{k}\leq \mathrm{b}-\mathrm{T}\mathrm{\bar{x}}\qquad&&\forall k\in[K].
\end{alignat}
\end{subequations}
We note that since $\mathrm{x}\in \mathcal{X}\cap \mathcal{V}$, $\mathtt{MP}$ enjoys \emph{relatively complete recourse}, \ie it has feasible solutions for all $\mathrm{\bar{x}} \in \mathcal{X}\cap \mathcal{V}$, and $\upxi \in \Xi$. Without this property, the $\mathtt{SP}$ might be infeasible for some $\bar{x}$, thus requiring feasibility cuts to be generated.  

The basic idea of the classical Benders algorithm is to approximate the function $\nu(\mathrm{x})$ using hyper-planes (referred to as \emph{optimality cuts}) generated by solving the dual $\mathtt{SP}$ for fixed values of \(\mathrm{x}\). However, since $\mathtt{SP}$ has binary decision variables, it is not possible to use the duality theory to generate cuts. Assuming that we have an oracle to solve $\mathtt{SP}$, in any iteration $r$, the $r$-th feasible solution $\mathrm{x}^{r}$ is used to define the sets 
$\mathcal{S}_{r}:=\{i\in[n]:x^r_i=1\}$
and to evaluate its corresponding worst-case second-stage objective function value \(\theta_{r}\). We use this solution and value to generate the valid combinatorial cut, first proposed by \cite{laporte1993integer}, \[\theta \geq (\theta_{r}-L_{r})(\sum_{i\in \mathcal{S}_{r}}{x_{i}}-\sum_{\notin  \mathcal{S}_{r}}{x_{i}})-(\theta_{r}-L_{r})(|\mathcal{S}_{r}|-1)+L_{r},\] where \(|\mathcal{S}_{r}|\) is the cardinality of \(\mathcal{S}_{r}\), \(L\) is a lower bound on the optimal value of the SP and \(\theta\) is a $\mathtt{MP}$ decision variable that defines the epigraph of $\nu(\mathrm{x})$. Hence, $\mathtt{MP}$ can be written as follows:
\begin{subequations}\label{MP}
\begin{alignat}{3}
&\;\min_{\mathrm{x} \in \mathcal{X},\theta}\;&& {\mathrm{c}'\mathrm{x}+\theta}\label{conMP1}\\
&\;\text{s.t.}\;&&\theta \geq (\theta_{r}-L_{r})(\sum_{i\in \mathcal{S}_{r}}{x_{i}}-\sum_{\notin  \mathcal{S}_{r}}{x_{i}})-(\theta_{r}-L_{r})(|\mathcal{S}_{r}|-1)+L_{r}\;\;\;\;\forall r\label{conMP2},
\end{alignat}
\end{subequations}
The Logic-based Benders decomposition algorithm is summarized as follows:

\begin{algorithm}[H]
\SetAlgoLined
\DontPrintSemicolon
 Initiate with an arbitrary feasible \(\mathrm{\bar{x}}\), set \(UB=\infty\), \(LB=-\infty\), \(r=1\)\; 
 
 \While{$UB-LB \ge \epsilon$}{

Solve $\mathtt{SP}$ \eqref{sub-problem} with \(\mathrm{\bar{x}}\) and find \(\theta_{r}\)\; 

Find a lower bound \(L_{r}\) for $\mathtt{SP}$ \eqref{sub-problem} as will be explained later\;

Generate the optimality cut: \(\theta \geq (\theta_{r}-L_{r})(\sum_{i\in \mathcal{S}_{r}}{x_{i}}-\sum_{i\notin \mathcal{S}_{r}}{x_{i}})-(\theta_{r}-L_{r})(|\mathcal{S}_{r}|-1|)+L_{r}\)\;

Update the upper bound: \(UB=\min{(UB,\mathrm{c}'\mathrm{\bar{x}}+\theta_{r})}\)\;

Solve $\mathtt{MP}$ \eqref{MP} with the new optimality cut added\;

Set \(LB\) equal to the optimal value of $\mathtt{MP}$ \eqref{MP}\;

Extract the optimal partial solution \(\mathrm{x}^{*}\) and use it as \(\mathrm{\bar{x}}\) in the next iteration and set \(r=r+1\)\;
 
 }
 
\textbf{Return:}{Declare the pair \(\left(\mathrm{x}^{*},\mathrm{y}^{{k}^{*}}\right)\) as the optimal solution}
\caption{Logic-based Benders decomposition algorithm for solving $K$-adaptability problem}\label{algo:lbbda}
\end{algorithm}

We note that $\mathtt{SP}$ \eqref{sub-problem} is a MMMCRO problem. In the next section, we propose a novel approach to solve it and also to obtain the lower bound $L_{r}$.

\subsection{Double-oracle for Solving SP}\label{doracle}
We define the partial sets $\mathcal{Y}'\subseteq \mathcal{Y}$ with $|\mathcal{Y}|\ge K$ and $\Xi'\subset \Xi$ and use $j\in [J]$ ($[J']$) and $h\in [H]$ ($[H']$) to index the elements of $\mathcal{Y}$ ($\mathcal{Y}'$) and the vertices of $\Xi$ ($\Xi'$), also referred to as ``scenarios", respectively. $\mathtt{SP}$ \eqref{sub-problem} can be reformulated over these partial sets as the $p$-center problem $\mathtt{P}(\mathcal{Y}',\Xi')$ \citep{arslan2020min}: 
\begin{subequations}\label{P5}
\begin{alignat}{3}
&\min_{\{z_{j}\}_{j\in[J']}, \{v_{jh}\}_{j\in[J'],h\in[H']}, w}\quad&&w\\
&\quad\qquad\qquad\text{s.t.}\quad&&w\ge \sum_{j\in [J']}\upxi_{h}'\mathrm{Q}\mathrm{y}_{j}v_{jh}\quad&&\forall h\in [H']\label{con11}\\
&&&\sum_{j\in [J']}v_{jh}=1\quad&&\forall h\in [H']\label{con21}\\
&&&\sum_{j\in [J']}z_j= K\label{con31}\\
&&&v_{jh}\le z_j\quad&&\forall j\in [J'], \forall h\in [H']\label{con41}\\
&&&v_{jh}, z_j\in \{0,1\}\quad&&\forall j\in [J'], \forall h\in [H'].
\end{alignat}
\end{subequations}
The binary variable $z_j$ takes value 1 if the feasible solution $\mathrm{y}_j$ is selected to be among the $K$ ``prepared" solutions, and $v_{jh}$ takes value 1 if scenario $\upxi_h$ is assigned to solution $\mathrm{y}_j$, and 0 otherwise. Constraint \eqref{con11} finds the scenario-solution pair with the worst cost among all assignments. Constraint \eqref{con21} ensures that each scenario is assigned to exactly one solution, whereas \eqref{con31} and \eqref{con41}, respectively, stipulate that $K$ solutions are selected and that scenarios can be assigned to selected solution only.

To solve \eqref{P5}, the following algorithm is proposed:

\begin{enumerate}
	\item Solve the problem $\mathtt{P}(\mathcal{Y}',\Xi)$, \ie the problem with the subset $\mathcal{Y}'$ of all recourse solutions generated so far (carried forward from Step 2 in the previous iteration) and all scenarios in $\Xi$ to obtain an upper bound $UB$. To solve this problem, we begin with a subset $\Xi'$ of scenarios and perform the following steps.
	\begin{enumerate}
		\item Solve the problem $\mathtt{P}(\mathcal{Y}',\Xi')$ (\ie Problem \eqref{P5}) to find $\mathrm{z}^*$, $\mathrm{v}^{*}$ and $w^*$. Identify the optimal $K$-subset of recourses as $\left\{\mathrm{y}^k\in \mathcal{Y}': z_k^*=1\right\}$
		\item Given the current optimal $K$-subset $\left\{\mathrm{y}^k\right\}_{k\in[K]}$ of recourses, try to find a scenario $\mathrm{\upxi}_{|H'|+1}\in \Xi$ that violates \eqref{con11} by solving the problem
		\begin{subequations}\label{senariogen}
		\begin{alignat}{3}
			&\max_{\upxi\in \Xi,\eta}\quad&&\eta\\
			&\;\;\;\text{s.t.}\quad&&\eta \le \upxi'\mathrm{Q}\mathrm{y}^{k}\quad&&k\in [K]	
		\end{alignat}
		\end{subequations}
		\item If $\eta^*>w^*$, add the new scenario to $\Xi'$ and repeat steps (a) and (b). Otherwise, stop and move to Step 2.
	\end{enumerate}

\item In this step, we find the optimal $K$-subset $\left\{{\mathrm{y}^k}^*\right\}_{k\in [K]}$ of recourses that minimizes the worst-case loss for the discrete scenario set $\Xi'$ by solving the problem $\mathtt{P}(\mathcal{Y},\Xi')$:
\begin{subequations}\label{P7}
\begin{alignat}{4}
&\min_{\{\mathrm{y}^k\}_{k\in[K]},\gamma,\{u_{kh}\}_{k\in[K],h\in[H']}}\quad&&\gamma \label{con61}\\
&\qquad\qquad\quad\text{s.t.}\quad&&\upxi_{h}'\mathrm{Q}\mathrm{y}^{k}\le \gamma+M(1-u_{kh})\quad&&\forall k\in[K],\,\forall h\in [H'] \label{con62}\\
&&&\sum_{k\in [K]}u_{kh}=1\quad&&\forall h\in [H'] \label{con63}\\
&&&\mathrm{y}^k\in \mathcal{Y}\quad&&\forall k\in [K] \label{con64}\\
&&&\mathrm{W}\mathrm{y}^{k}\leq \mathrm{b}-\mathrm{T}\mathrm{\bar{x}}\quad&&\forall k\in [K] \label{con65}\\
&&&u_{kh}\in  \{0,1\}\quad&&\forall k\in[K],\,\forall h\in [H'] \label{con66}.
\end{alignat}
\end{subequations}
In problem \eqref{P7}, $u_{kh}$ is a binary assignment variable which takes value 1 if scenario $h$ is assigned to recourse $k$. We update the solution pool as \(\mathcal{Y}' \leftarrow \mathcal{Y}' \cup \left\{{\mathrm{y}^k}^*\right\}_{k\in [K]}\), where \(\left\{{\mathrm{y}^k}^*\right\}_{k\in [K]}\) is the partial optimal solution of problem \eqref{P7}. Moreover we set $LB=\gamma$.

\item Iterate between steps (1) and (2) until $UB-LB<\varepsilon$. Declare the incumbent $\left\{{\mathrm{y}^k}^*\right\}_{k\in[K]}$ as the optimal solution.
\end{enumerate}

\begin{algorithm}[H]
\SetAlgoLined
\DontPrintSemicolon
 initialization: $\mathrm{y}' , \Xi', LB=-\infty, UB=+\infty$\; 
 
 \While{$UB-LB \ge \epsilon$}{

 \While{Scenario-added=\textbf{true}}{
 
 Compute $w^*$, $\mathrm{z}^*$, $\mathrm{v}^*$, and ${\mathrm{y}^{k}}^*$ by solving \eqref{P5}\;
 
 Compute $\xi_{|H|+1} \in \Xi$, and $\eta^*$ by solving \eqref{senariogen}\;
 
 \eIf{$\eta^*>w^*$}
 {
$\Xi' \leftarrow \Xi' \cup \left\{\xi_{|H|+1}\right\}_{k\in [K]}$\;

Scenario-added=\textbf{true}

 }{$\Xi' \leftarrow \Xi'$\;
 
 Scenario-added=\textbf{false}\;
 }   
 }

\textbf{Return:} $\Xi'$, $UB=w^{*}$

Compute $\{{y^k}^*\}_{k\in[K]}$, $\gamma^*$ by solving \eqref{P7}\; 

$\mathcal{Y}' \leftarrow \mathcal{Y}' \cup \left\{{\mathrm{y}^k}^*\right\}_{k\in [K]}$\;

$LB=\gamma^*$\; 
 }
\textbf{Return:} $\{{y^{k}}^*\}_{k\in[K]}$ 
\caption{The Double-Oracle algorithm for solving $\mathtt{SP}$ \eqref{sub-problem}}
\end{algorithm}

To generate an optimality cut, the logic-based Benders decomposition algorithm explained in the previous section requires a valid lower bound ($L$) on the optimal value of the second-stage problem. A valid lower bound should be \(L\leq \min\limits_{\mathrm{x}}\{\nu(\mathrm{x})|\;\mathrm{x}\in \mathcal{X}\}\). In every iteration of the proposed algorithm, we calculate a lower bound on the optimal value of $\mathtt{SP}$ for a fixed first-stage decision \(\mathrm{\bar{x}}\) by solving \eqref{P7}. However, a lower bound on the optimal value of $\mathtt{SP}$ for all \(\mathrm{x}\in \mathcal{X}\) is required, which can be obtained by solving the following problem: 
\begin{subequations}\label{Lower-bound}
\begin{alignat}{4}
&\min_{\{\mathrm{y}^k\}_{k\in[K]},\gamma,\{u_{kh}\}_{k\in[K],h\in[H']}, \mathrm{x}}\quad&&\gamma \label{conlb61}\\
&\;\qquad\qquad\quad\text{s.t.}\quad&&\upxi_{h}'\mathrm{Q}\mathrm{y}^{k}\le \gamma+M(1-u_{kh})\quad&&\forall k\in[K],\,\forall h\in [H'] \label{conlb62}\\
&&&\sum_{k\in [K]}u_{kh}=1\quad&&\forall h\in [H'] \label{conlb63}\\
&&&\mathrm{y}^k\in \mathcal{Y}\quad&&\forall k\in [K] \label{conlb64}\\
&&&\mathrm{W}\mathrm{y}^{k}\leq \mathrm{b}-\mathrm{T}\mathrm{x}\quad&&\forall k\in [K] \label{conlb65}\\
&&&\mathrm{x}\in \mathcal{X} \label{conlb67}\\
&&&u_{kh}\in  \{0,1\}\quad&&\forall k\in[K],\,\forall h\in [H'] \label{conlb66}.
\end{alignat}
\end{subequations}
Note that the lower bond ($L_{r}$) changes at each Benders iteration since we solve problem \eqref{Lower-bound} in each iteration by using an updated subset of scenarios $\Xi'$. We also suggest “warm-starting” the SP in every Benders iteration by re-using some of the $\mathrm{y}$ variables generated in previous iterations. Given a subset $\left\{\mathrm{y}_j\right\}_{j\in [J']}$ of recourse solutions, one can “filter” them using constraint \eqref{con65} and reuse the ones that satisfy this constraint for the new $\bar{\mathrm{x}}$ in problem \eqref{P5} right away. Likewise, the scenarios (vertices of $\Xi$) generated in an iteration can be re-used in subsequent iterations of the Benders algorithm vertices since they do not depend on $\bar{\mathrm{x}}$. Such warm-staring techniques can substantially improve the performance of the proposed algorithm, even though we have not use them in our numerical tests. Figure \ref{fig:alg} illustrates the proposed algorithm:

\begin{figure}[H]
    \centering
    \includegraphics[width=12cm]{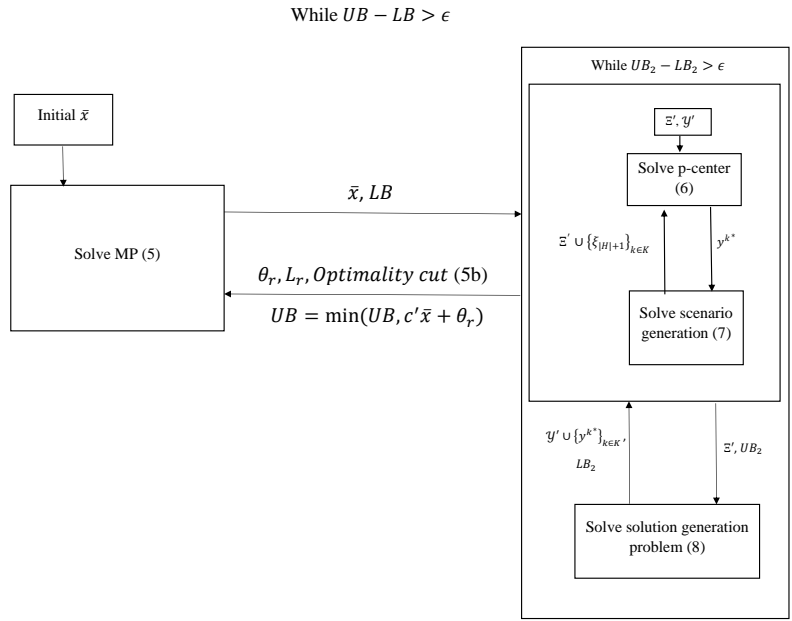}
     \caption{Double-oracle, Logic-based Benders decomposition algorithm}
     \label{fig:alg}
\end{figure}

\section{Finite Convergence}\label{sec3}

The proposed algorithm consists of three loops. The outer loop handles the first-stage decision variables by using logic-based Benders decomposition. The other two loops are used for solving $\mathtt{SP}$, which is a MMMRCO problem. The inner loop has two step: scenario generation and solution generation. We show that the outer and the scenario generation loops terminate in a finite number of iterations, and the solution generation loop leads to convergence of \(LB\) and \(UB\). These lemmas do not exploit the linearity of the cost and recourse constraint functions in problem \eqref{K-adaptability}.
 
\begin{lemma}
Given that \(\mathcal{X}\) is bounded, the number of optimality cuts \eqref{conMP2} generated and, thus, the number of iterations of the outer loop are finite.
\end{lemma}

\begin{proof}
Since \(\mathcal{X}\) is bounded and \(\mathrm{x}\) are binary, there is a finite number of feasible solutions $\bar{\mathrm{x}}$ for the first-stage problem \eqref{MP}. Each first-stage solution corresponds to a single optimality cut \eqref{conMP2}, generated by solving the the $\mathtt{SP}$ to obtain $\mathrm{x}^r$, \(\theta_{r}\) and \(L\). Hence, the number of optimality cuts is finite, and so is the number of outer loop iterations.
\end{proof}

\begin{lemma}
The maximum number of scenarios $\left\{\upxi_h\right\}_{h\in H}$ that can be generated through the double oracle algorithm is finite, and equal to $\frac{|J|!}{K!(|J|-K)!}$. 
\end{lemma}

\begin{proof}
Given that \(\mathcal{Y}\) is a bounded discrete set, its elements (feasible solutions) are finite and thus can be enumerated. The $p$-center problem \eqref{P5} selects $K$ solutions and assigns those selected solutions to scenarios. Moreover, those $K$ selected solutions are used to generate worst-case scenarios through problem \eqref{senariogen}; and if the generated scenario violates constraint \eqref{con11}, it is added to the set of scenarios. However, if all possible worst-case scenarios of all feasible \((\mathrm{\bar{x}},\mathrm{y}^k)\) assignment in our subset of scenarios are available, then no generated scenario from problem \eqref{senariogen} can violate \eqref{con11} and the iteration between \eqref{P5} and \eqref{senariogen} will terminate. If all combination of $K$ out of $|\mathcal{Y}|$ are selected and used in problem \eqref{senariogen}, all possible worst-case scenarios for all feasible \((\mathrm{\bar{x}},\mathrm{y}^k)\) are generated. Since \(\mathcal{Y}\) is a bounded discrete set, there are exactly \(\frac{|J|!}{k!(|J|-k)!}\) possible ways to select $K$ elements from the set \(\mathcal{Y}\). Consequently, at most there are \(\binom{|J|}{k} = \frac{|J|!}{k!(|J|-k)!}\) possible solutions for \eqref{P5}. Hence, in the worst-case situation, with the finite number of iterations between the $p$-center problem, \eqref{P5}, and the scenario generation problem, \eqref{senariogen}, all possible worst-case scenarios are generated. therefore, no new scenario can violate \eqref{con11} and then, this step terminates.    
\end{proof}

\begin{lemma}
The upper and lowed bounds, obtained respectively by solving problems $P(\mathcal{Y}',\Xi)$ and $P(\mathcal{Y},\Xi')$, will converge.
\end{lemma}

\begin{proof}
The upper bound is achieved by solving $\mathtt{P}(\mathcal{Y}',\Xi)$. In order to get the optimal assignment, problem \eqref{P5} with $\mathcal{Y}'$ must be solved. However, the optimal value of its objective function, which is the upper bound, can be achieved by solving \(UB=\max_{h \in [H]} \min_{j \in [J']}\upxi'_{h}\mathrm{Q}\mathrm{y}_{j}\), the proof of the reformulation is provided in the appendix A. On the other hand, By using the fixed set of scenarios, a set of $K$ solutions are generated by solving \eqref{P7}, which can be reformulated as \(LB=\max_{h \in [H']} \min_{\mathrm{y}^{k}\in \mathcal{Y}}\upxi'_h\mathrm{Q}\mathrm{y}^{k}\), the detail of this reformulation is provided in the appendix A. There are two cases.

First: The optimal solutions \({\mathrm{y}^k}^{*}\) already exist in the subset of solutions $\mathcal{Y}'$, in this case, the optimal pair of solutions and scenarios will be found by solving \eqref{P5} and its objective value is equivalent to \(UB=\max_{h \in [H]} \min_{j \in [J']}\upxi'_{h}\mathrm{Q}\mathrm{y}_{j}\). On the other hand, the optimization problem related to lower bound, \eqref{P7}, will generate \({\mathrm{y}^k}^{*}\) as optimal solution since the scenarios are fixed for both the lower and upper bound problems. In this case $LB=UB$ and the algorithm will terminate. 

Second: The optimal solutions \({\mathrm{y}^k}^{*}\) are not in the subset of solutions $\mathcal{Y}'$. In this case, by using fixed scenarios and solving \eqref{P7} a set of $K$ solutions are generated that includes \({\mathrm{y}^k}^{*}\). This set of solutions are added to the subset of solutions in the $p$-center problem \eqref{P5}. Consequently, the optimal pair of solutions and scenarios will be in the \(\mathcal{Y}'\), and \(\Xi'\), respectively. Hence, by solving the $p$-center problem \eqref{P5}, the optimal pair of solutions and scenarios will be selected. Consequently, $LB=UB$ and the algorithm will terminate.
\end{proof}

Based on above mentioned lemmas, we can conclude that all three loops in proposed algorithm will terminate in finite number of iterations. Consequently, the proposed algorithm will converge in finite number of iterations. 

\section{Extensions}\label{sec4}

So far, we focused on the $K$-adaptability problem with linear objectives and constraints, binary decision variables in the first- and second-stage that can be extend to integer second-stage decision variables, and with objective uncertainty only in the second-stage. In this section, we show how the proposed algorithm can be extended to more general cases beyond the basic setting outlined earlier.

\subsection{Second-stage Constraint Uncertainty}

Similar to the algorithm proposed by \cite{arslan2020min} to solve the MMMRCO problem, our approach can be extended to the case when uncertainty affects both the objective function and constraints of the recourse problem. The extended problem can be formulated as follows:
\begin{subequations}\label{K-adaptability2}
\begin{alignat}{3}
&&\;\min_{\mathrm{x} \in \mathcal{X}, \{\mathrm{y}^{k}\}_{k\in [K]}}&\mathrm{c}'\mathrm{x}+\sup_{\upxi \in \Xi}\min_{k\in [K]}\{\upxi'\mathrm{Q}\mathrm{y}^{k}:\;\mathrm{T}\mathrm{x}+\mathrm{W}(\upxi)\mathrm{y}^{k}\leq \mathrm{b}\}\\
&&\;\text{s.t.}\;\;\;\;\;\; &\mathrm{y}^{k}\in \mathcal{Y},\;k\in [K],
\end{alignat}
\end{subequations}
where \(\mathrm{W}(\upxi)\) is an affine mapping of uncertain parameters. Note the dependency of the $\mathtt{SP}$ constraints on \(\upxi\). To solve this problem, we use the same iterative algorithm explained earlier but with a modified $p$-center problem $\mathtt{P}'(\mathcal{Y}',\Xi')$ by adding \eqref{conP-center26}, as follows:
\begin{subequations}\label{P-center2}
\begin{alignat}{4}
&\min_{\{z_{j}\}_{j\in [J']}, \{v_{jh}\}_{j\in [J'],h\in [H']}, w}\quad&&w\label{conP-center21}\\
&\qquad\qquad\;\;\;\text{s.t.}\quad&& w\ge \sum_{j\in [J']}\upxi'_{h}\mathrm{Q}\mathrm{y}_{j}v_{jh}\quad&&\forall h\in [H']\label{conP-center22}\\
&&&\sum_{j\in [J']}v_{jh}=1\quad&&\forall h\in [H']\label{conP-center23}\\
&&&\sum_{j\in [J']}z_j= K\label{conP-center24}\\
&&&v_{jh}\le z_j\quad&&\forall j\in [J'], \forall h\in [H']\label{conP-center25}\\
&&&v_{jh}=0\quad&&\forall j\in [J'],\;h\in [H']: \exists s\in [S]\nonumber\\
&&&&&\text{s.t.}\; \mathrm{e'}_{s}(\mathrm{W}(\upxi_{h})\mathrm{y}_{j}-\mathrm{b}+\mathrm{T}\mathrm{\bar{x}})> 0\label{conP-center26} \\
&&&v_{jh}, z_j\in \{0,1\}\quad&&\forall j\in [J'], \forall h\in [H']\label{conP-center27}.
\end{alignat}
\end{subequations}
In this formulation, there are $s$ constraints with uncertain parameters, \ie \(\mathrm{b} \in \mathbb{R}^{s}\), and \(\mathrm{e}_{s}\) denotes the $s$-th column of the identity matrix \(\mathrm{I}_{s}\). Constraint \eqref{conP-center26} prevents infeasible assignment of solution-scenario pairs. Moreover, problem \eqref{senariogen} is modified by adding \eqref{modification} as follows:
\begin{subequations}\label{senariogen2}
\begin{alignat}{3}
&\max_{\xi\in \Xi,\eta, \lambda\in\{0,1\}^K}\quad&&\eta\\
&\qquad\;\;\text{s.t.}\quad&&\eta \le \upxi'\mathrm{Q}\mathrm{y}^{k}+M\lambda^{k}\quad&&k\in [K]\\
&\;\;\; \quad&&\mathrm{W}(\upxi)\mathrm{y}^{k}\geq \mathrm{b}-\mathrm{T}\mathrm{\bar{x}}-M(1-\lambda^{k})+\epsilon\quad&&k\in [K]\label{modification},
\end{alignat}
\end{subequations}
where \(\lambda^{k}=1\) if the scenario $\xi$ is such that $\mathrm{y}^k$ is infeasible for any of the \(s\) uncertain constraints, thus enforces that we do not consider \(\upxi'\mathrm{Q}\mathrm{y}^{k}\) for calculating the upper bound on $\eta$.
In the case where \(\mathrm{W}(\upxi)\mathrm{y}^{k}= \mathrm{b}-\mathrm{T}\mathrm{\bar{x}}\), we would want \(\lambda^k=0\). Consequently, the small \(\epsilon\) on the RHS of \eqref{modification} prevents \( \lambda^k=1\) for equality case.
Moreover, the solution generation problem \eqref{P7} is modified by considering constraint \eqref{consg65} with uncertain parameters.
\begin{subequations}\label{solution-generation}
\begin{alignat}{4}
&\min_{\{\mathrm{y}^k\}_{k\in[K]},\gamma,\{u_{kh}\}_{k\in[K],h\in[H']}}\quad&&\gamma \label{consg61}\\
&\qquad\;\;\;\text{s.t.}\quad&&\upxi_{h}'\mathrm{Q}\mathrm{y}^{k}\le \gamma+M(1-u_{kh})\quad&&\forall k\in[K],\,\forall h\in [H'] \label{consg62}\\
&&&\sum_{k\in [K]}u_{kh}=1\quad&&\forall h\in [H'] \label{consg63}\\
&&&\mathrm{y}^k\in \mathcal{Y}\quad&&\forall k\in [K] \label{consg64}\\
&&&\mathrm{W}(\upxi'_{h})\mathrm{y}^{k}\leq \mathrm{b}-\mathrm{T}\mathrm{\bar{x}}+M(1-u_{kh})\quad&&\forall k\in [K],\forall h\in [H'] \label{consg65}\\
&&&u_{kh}\in  \{0,1\}\quad&&\forall k\in[K],\,\forall h\in [H'] \label{consg66}.
\end{alignat}
\end{subequations}

One should note that, even though $\mathcal{V}$ can be modified to ensure that $\mathtt{MP}$ satisfies relatively complete recourse, it might happen that for some $\bar{\mathrm{x}}\in\mathcal{X}\cap \mathcal{V}$, there exist no set of $K$ solutions $\{\mathrm{y}^k\}_{k=1}^K$ that ensure that some $\mathrm{y}^k$ is always feasible under all $\upxi\in\Xi$. For this reason, at the $r$-th iteration of the the logic-based Benders decomposition algorithm, $\mathtt{SP}$ might become infeasible for $\mathrm{x}_i^r$, which is identified when problem \eqref{solution-generation} becomes infeasible for some $\Xi'$. At this point, Algorithm \ref{algo:lbbda} should be modified to  returns to the $\mathtt{MP}$ a feasibility cut of the form:
\[\sum_{i\in \mathcal{S}_{r}}{\mathrm{x}_{i}}-\sum_{\notin  \mathcal{S}_{r}}{\mathrm{x}_{i}}\leq |\mathcal{S}_{r}|-1,\]
in order to discard $\mathrm{x}_i^r$ from the set of feasible candidates, instead of producing an optimality cut of the form \eqref{conMP2}.

\subsection{First-stage Integer Decision Variables}
The outer loop in the proposed algorithm depends on the first-stage variables being binary (\ie \(\mathrm{x} \in \{0,1\}^{n}\)) to generate logic-based Benders cuts of the type \eqref{conMP2}. However, if the first-stage integer variables are not binary, but rather general integer, \ie \(\mathcal{X} \in \mathbb{Z}^{n}_{+}\), one can simply apply the transformation \(x_i=\sum_{p_i=0}^{P_i}{2^pu_{ip_i}},\; i=1,\dots,n\), where \(u_{ip_i}\in \{0,1\}\), and \(P_i\) depends on upper bound of $x_i$. Clearly, this generalization comes at the expense of increasing the number of variables in the first-stage problem, thus it might be efficient only for small values of $P$. It should be noted that the basic algorithm described in Section \ref{sec2} can handle general integer recourse decisions since none of the algorithm steps depends on $\mathrm{y^k}$ being binary.


\subsection{First-stage Objective Uncertainty}

Even though \cite{bertsimas2010finite} defined $K$-adaptability such that the first-stage objective is deterministic, we extend our algorithm to the case when there is first-stage objective uncertainty, similar to \cite{hanasusanto2015k}. A practical example of the $K$-adaptability problem with uncertainty in the first-stage is the multi-period portfolio selection problem where decisions about the allocation of budget among assets have to be made at the beginning of the investment horizon, thus are first-stage decisions. Indeed, asset returns are uncertain even at the outset, and the initial capital allocation decision cannot be postpone until this uncertainty is revealed. Hence, the here-and-now decisions are directly affected by the uncertain parameters. In this section, we differentiate between two cases of first-stage objective uncertainty: when it is independent from the second-stage uncertainty, and when some (or all) uncertain parameters affect both stages, \ie dependent uncertainty. We show how the proposed approach is tailored for each case.

\subsubsection{Independent Uncertainty}\label{independent}

Let us consider the following $K$-adaptability problem:   
\begin{subequations}\label{K-adaptability first-stage uncertainty}
\begin{alignat}{3}
&\min_{\mathrm{x} \in \mathcal{X}}\max_{\mathrm{\upxi} \in \mathrm{\Xi}, \upomega\in \Omega}\min_{k\in K}\{\upomega'\mathrm{C}\mathrm{x}+\mathrm{\upxi}'\mathrm{Q}\mathrm{y}^{k}:\;\mathrm{T}\mathrm{x}+\mathrm{W}\mathrm{y}^{k}\leq \mathrm{b}\}\\
&\;\text{s.t.}\;\; \; \mathrm{y}^{k}\in \mathcal{Y},\;k\in [K],
\end{alignat}
\end{subequations}
where \(\upomega \in \mathbb{R}^{c}\), \(C \in \mathbb{R}^{c\times n}\), and \(\Omega\) is a compact and convex uncertainty set. Other variables and parameters are the same as those used in formulation \eqref{K-adaptability}. We assume that \(\upomega\) and \(\upxi\) are disjoint sets of uncertain parameters. In this case, $\mathtt{MP}$ \eqref{MP} is modified as follows:

\begin{subequations}\label{MP first-stage uncertainty}
\begin{alignat}{3}
&\min_{\mathrm{x} \in \mathcal{X},\; \theta}\max_{\upomega \in \Omega} {\upomega'\mathrm{C}\mathrm{x}+\theta}\label{conMPuf1}\\
&\;\text{s.t.}\; \theta \geq (\theta_{r}-L)(\sum_{i\in \mathcal{S}_{r}}{x_{i}}-\sum_{\notin  \mathcal{S}_{r}}{x_{i}})-(\theta_{r}-L)(|\mathcal{S}_{r}|-1)+L\label{conMPuf2},
\end{alignat}
\end{subequations}
which is a static RO problem that can be tractably reformulated by applying convex duality on the inner maximization. The SP does not change, and we can apply the double-oracle algorithm described in Section \ref{doracle} to solve it. 

\subsubsection{Dependent Uncertainty}

Next, we consider the $K$-adaptability problem variant addressed by \citep{hanasusanto2015k}, where the same uncertain parameters affect both stages, formulated as follows:  
\begin{subequations}\label{K-adaptability first-stage uncertainty2}
\begin{alignat}{3}
&&\min_{\mathrm{x} \in \mathcal{X},\{y^k\}_{k\in[K]}}& \max_{\mathrm{\upxi} \in \mathrm{\Xi}}\min_{k\in K}\upxi'\mathrm{C}\mathrm{x}+\mathrm{\upxi}'\mathrm{Q}\mathrm{y}^{k}\\
&&\;\text{s.t.}\;\;\;\;\;\; \;& \mathrm{y}^{k}\in \mathcal{Y},\;\mathrm{T}x+\mathrm{W}\mathrm{y}^k\leq \mathrm{b},\;k\in [K],
\end{alignat}
\end{subequations}
where \(\mathrm{C} \in \mathbb{R}^{q\times n}\) and the rest of parameters and variables are the same as in \eqref{K-adaptability}. This case can be reformulate as second stage uncertainty. 

Let \(\bar{\mathrm{Q}}=[\mathrm{C},\mathrm{Q}] \in \mathbb{R}^{q\times (n+m)}\), $\bar{\mathrm{T}}=\begin{bmatrix}
           \mathrm{I} \\
           \mathrm{I} \\
           \mathrm{T}
         \end{bmatrix}$, $\bar{\mathrm{W}}=\begin{bmatrix}
           -\mathrm{I} & 0 \\
           \mathrm{I} & 0 \\
           0 & \mathrm{W}
         \end{bmatrix}$, and  \(\bar{b}=\begin{bmatrix}
           0 \\
           0 \\
           \mathrm{b}
         \end{bmatrix}\). We have that problem \eqref{K-adaptability first-stage uncertainty2} can be equivalently formulated as:
\begin{subequations}\label{K-adaptability first-stage uncertainty3}
\begin{alignat}{3}
&&\min_{\mathrm{x} \in \mathcal{X},\{\bar{\mathrm{y}}^k\}_{k\in[K]}}&\max_{\mathrm{\upxi} \in \mathrm{\Xi}}\min_{k\in K}\upxi'\bar{\mathrm{Q}}\bar{\mathrm{y}}^{k}\\
&&\text{s.t.}\;\;\;\;\;\;&\bar{\mathrm{y}}\in\{0,1\}^n\times \mathcal{Y},\;\bar{\mathrm{T}}x+\bar{\mathrm{W}}\bar{\mathrm{y}}^k\leq \bar{\mathrm{b}},\,k\in[K].
\end{alignat}
\end{subequations}

\subsection{Nonlinear Objective and Constraint Functions}\label{nonlinear}

Our algorithm can handle nonlinear $K$-adaptability problems of the form:  
\begin{subequations}\label{nonlinear K-adaptability}
\begin{alignat}{3}
&\min_{\mathrm{x} \in \mathcal{X}}f(\mathrm{x})+\max_{\upxi \in \mathrm{\Xi}}\min_{k\in K}\{g(\upxi,\mathrm{y}^{k}):\;h(\upxi, \mathrm{x}, \mathrm{y}^{k})\leq \mathrm{b}\}\\
&\;\text{s.t.}\;\; \; \mathrm{y}^{k}\in \mathcal{Y},\;k\in [K],
\end{alignat}
\end{subequations}
where $f:\mathcal{X}\mapsto \mathbb{R}$ is convex in \(\mathrm{x}\),  $g:\Xi\times\mathcal{Y}\mapsto \mathbb{R}$ is affine in $\upxi$ and convex in $\mathrm{y}$, and  $h:\Xi\times\mathcal{X}\times\mathcal{Y}\mapsto \mathbb{R}$ is affine in \(\upxi\) and jointly convex in \(\mathrm{x}\) and \(\mathrm{y}\). In this case, $\mathtt{MP}$ \eqref{MP} is modified by using \(f(\mathrm{x})\) instead of \(\mathrm{c}'\mathrm{x}\). Moreover, the $p$-center problem \eqref{P-center2} is modified by replacing \(\upxi'_{h}\mathrm{Q}y_{j}v_{jh}\) and \(\upxi'\mathrm{W}y+\mathrm{T}\mathrm{\bar{x}}\) with \(g(\xi_{h},y_{j})\) and \(h(\xi_{h},y_{j},\mathrm{\bar{x}})\), respectively. In the $p$-center problem, \(v_{jh}\) and \(z_j\) are decision variables while \(\mathrm{y}_{j}\), \(\upxi_{h}\), and \(\mathrm{\bar{x}}\) are constant. Consequently, regardless of the type of functions \(g(\cdot,\cdot)\) and \(h(\cdot,\cdot,\cdot)\), the $p$-center problem finds the optimal $K$ solutions and assigns them to scenarios. Moreover, the scenario generation problem \eqref{senariogen} is rewritten as follows: 
\begin{subequations}\label{senariogen3}
\begin{alignat}{3}
&\max_{\xi\in \Xi,\eta,\lambda\in\{0,1\}^K}\quad&&\eta\\
&\qquad\;\;\text{s.t.}\quad&&\eta \le g(\xi,\mathrm{y}^{k})+M\lambda^{k}\quad&&k\in [K]\\
&\;\;\; \quad&& h(\xi,\mathrm{y}^{k},\bar{\mathrm{x}})\geq \mathrm{b}-M(1-\lambda^{k})+\epsilon\quad&&k\in [K]\label{senariogen3const2},
\end{alignat}
\end{subequations}
which again takes the form of a mixed integer LP given our assumption that $g(\cdot,\cdot)$ and $h(\cdot,\cdot,\cdot)$ be affine in $\upxi$. 
Finally, the solution generation problem \eqref{P7} is changed by replacing \(\upxi'_{h}\mathrm{Q}\mathrm{y}^{k}\) and \(\mathrm{W}\mathrm{y}^{k}+\mathrm{T}\mathrm{\bar{x}}\) with \(g(\upxi_{h},\mathrm{y}^{k})\) and \(h(\upxi_{h},\mathrm{y}^{k},\bar{\mathrm{x}})\), respectively.  

\section{Numerical Results}\label{sec5}

This section presents and discusses the numerical results obtained by implementing the proposed algorithm to solve three problems: the shortest path problem, the knapsack problem, and a generic $K$-adaptability problem. These results are compared to those obtained from state-of-the-art algorithms proposed in (MILP reformulation of \cite{hanasusanto2015k}, the iterative algorithm (IA) of \cite{chassein2019faster}, the row-and-column generation (RCG) algorithm proposed by \cite{goerigk2020min}, the scenario generation (SG) approach developed by \cite{arslan2020min}, and Branch-and-Bound (BB) method proposed by \cite{subramanyam2020k}). The proposed algorithm, Double-Oracle (DO), MILP, IA, and RCG were coded on Python 3.10.4 using Jupyter, SG is available \href{https://github.com/mjposs/min-max-min}{here}, and BB is also available \href{https://github.com/AnirudhSubramanyam/KAdaptabilitySolver/blob/v1.0/README.md}{here}). The subproblems were solved using CPLEX called through CPLEX-CMD on a Linux laptop with an 8th generation Intel Core i7 7700 processor and 16 GB RAM. The time limit and the relative optimality gap were set, respectively, to two hours (7200 seconds) and \(5\%\).

\subsection{Shortest Path Problem}

The first problem used to evaluate the proposed algorithm is the adaptive shortest path
problem, previously studied in \cite{hanasusanto2015k}, \cite{chassein2019faster}. This problem aims to select a subset of network arcs with the least total cost to form a path from a source $s$ to a destination $t$ when arc costs are uncertain. In the $K$-adaptability variant of the problem, $K$ paths are pre-formed and the shortest (least costly) among them is selected once the actual costs ate realized. We used test instances from \cite{arslan2020min}, available \href{https://github.com/mjposs/min-max-min}{here}.

Formally, the problem can be described as follows: A network \((\mathcal{V},\mathcal{A})\) has the cost of each arc \((i,j)\in \mathcal{A}\) characterized as \(c_{ij}=\bar{c}_{ij}+\xi_{ij}\hat{c}_{ij}\), where \(\bar{c}_{ij}\) is the nominal cost and \(\hat{c}_{ij}\) is the maximal deviation. The primary uncertain parameter \(\upxi\) belongs to the budgeted uncertainty set \(\Xi=\{\upxi \in [0,1]| \sum_{(i,j)\in A}\xi_{i,j}\leq \Gamma\}\), where \(\Gamma\) is an ``uncertainty budget" that controls the size of uncertainty set. With that, the problem is formulated as follows:

\begin{subequations}\label{P9}
\begin{alignat}{4}
&\min_{x_{ij} \in [0,1]^{n}}{\sum_{(i,j)\in \mathcal{A}}}{c_{ij}x_{ij}}\\
&\;\text{s.t.}\sum_{(i,j)\in \delta^{+}(i)}{x_{ij}}-\sum_{(i,j)\in \delta^{-}(i)}{x_{ij}}=b_{i},\;\forall i \in \mathcal{V},
\end{alignat}
\end{subequations}
where \(b_{s}=-1\), \(b_{t}=1\), and \(b_{i}=0\) for \(i\in \mathcal{V} /\ \{s,t\}\) and the sets \(\delta^{+}_{i}\) and \(\delta^{-}_{i}\)represent the forward and backward starts of node \(i \in \mathcal{V}\), respectively. 

We solve the problem in different sizes of \(|\mathcal{V}| \in \{20,25,40,50\}\). For each problem size, we considered \(k \in \{2,3,4,5,6\}\). Finally, each instance is solved based on different \(\Gamma \in \{3,6\}\). Ten randomly-generated instances were solved for each combination of $\nu$, $k$, and $\Gamma$. We compare the results of our algorithm with those based on the MILP reformulation of \cite{hanasusanto2015k}, the row-and-column generation (RCG) algorithm proposed by \cite{goerigk2020min}, the iterative algorithm (IA) of \cite{chassein2019faster} and the scenario generation (SG) approach developed by \cite{arslan2020min}.  

\begin{figure}[H]
    \centering
    \includegraphics[width=\textwidth,height=\textheight,keepaspectratio]{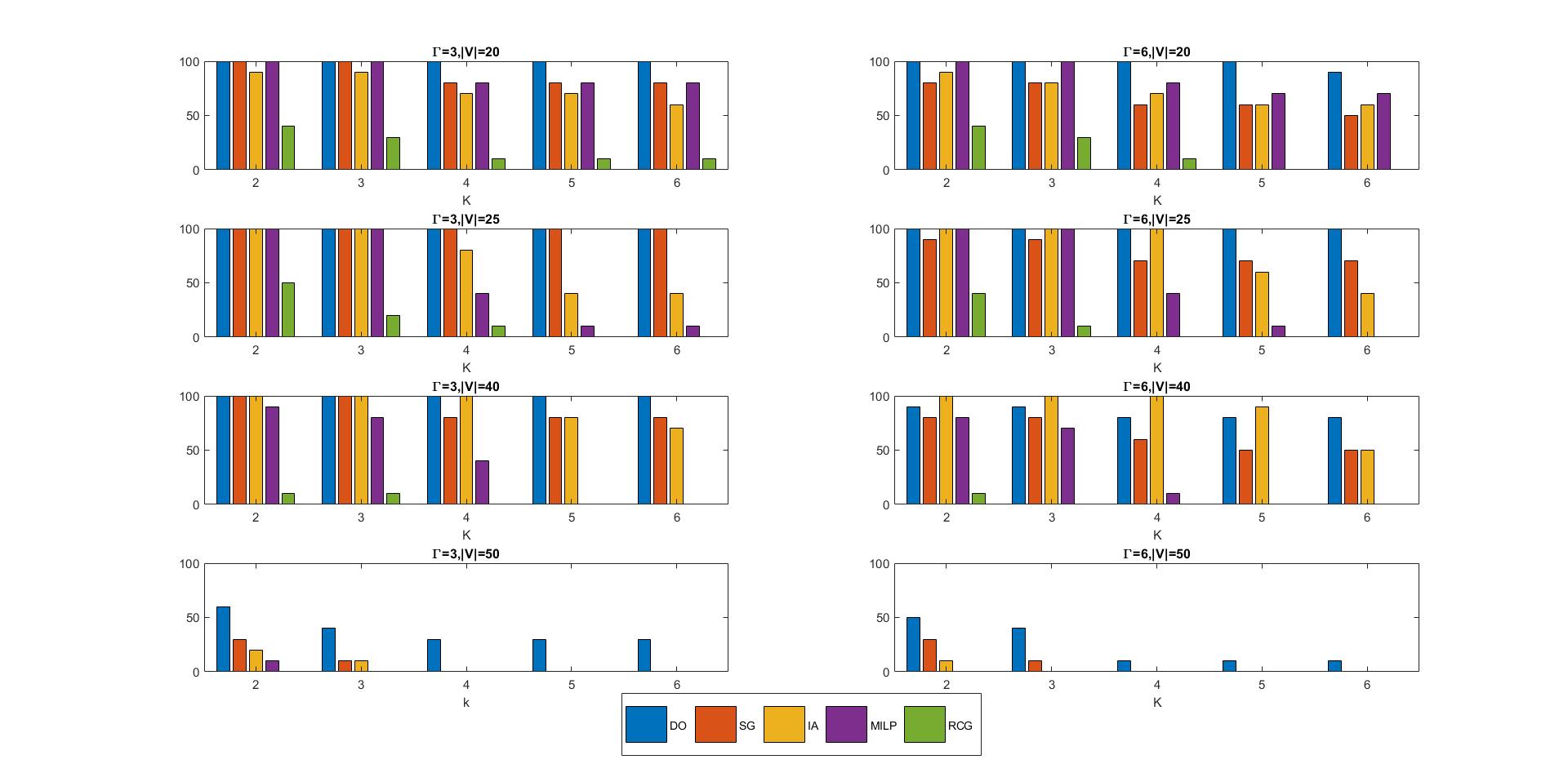}
     \caption{Percentage of solved instances of the shortest path problem}
     \label{fig:fig1}
\end{figure}

Figure \ref{fig:fig1} shows the percentage of instances solved by each algorithm within the cut-off time of two hours. Intuitively, as the problem size increases (in terms of both $|\mathcal{V}|$ and $K$), fewer instances are solved to optimality by all algorithms. Nevertheless, our proposed algorithm, labeled ``Double-oracle", shows better performance than all other algorithms. For example, our proposed algorithm solved all instances with \(\Gamma=3\) and \(k\in \{2,3,4,5\}\), while none of the other algorithms could solve all of these instances within the cut-off time. Moreover, the proposed algorithm solved $30-60\%$ of the instances with $|\mathcal{V}|=50$, $\Gamma=3$, and \(k\in \{2,3,4,5,6\}\), while the next best algorithm is the scenario generation method proposed by \cite{arslan2020min} that could not solve any instances with size $50$ and $K>3$. Details of the results can be found in Tables \ref{tab:multicol1} and \ref{tab:multicol2} in Appendix \ref{AppendixB}. 

The comparison results for \(\Gamma=6\), shown in Figure \ref{fig:fig1}, exhibit the same pattern. Our algorithm has a significant performance over the benchmark algorithms, except the iterative approach proposed by \cite{chassein2019faster}, when the uncertainty budget \(\Gamma\) is doubled. The iterative approach solved $100\%$ of instances with $|\mathcal{V}|=40$, $k\in \{2,3,4\}$, and $\Gamma=6$, while our algorithm solved $80\%$ to $90\%$ of instances with the same size. However, for $|\mathcal{V}|=50$, performance of the iterative algorithm deteriorates, as it could only solve $0-10\%$ of the instances, whereas the double-oracle algorithm solved $10-50\%$ of these instances. These results clearly show the performance advantage of the proposed algorithm over other algorithms proposed in the literature, especially for large-size problems. It is worthy to note that the iterative approach of \cite{chassein2019faster}, unlike ours, cannot guarantee global optimality since it uses a fixing heuristic to handle a bilinear term in each iteration.


\begin{figure}[H]
    \centering
    \includegraphics[width=\textwidth,height=\textheight,keepaspectratio]{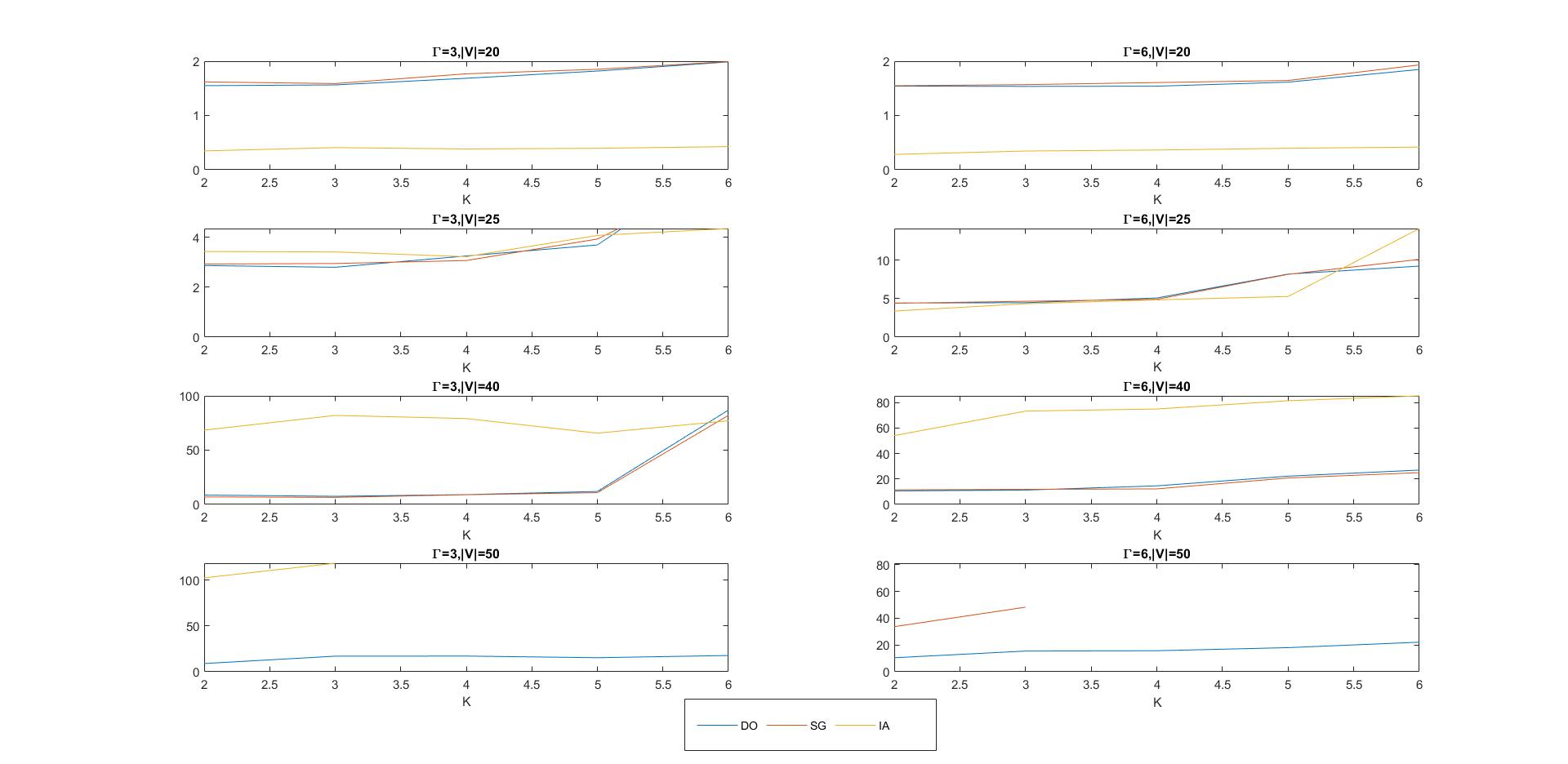}
     \caption{Average CPU times for the solved instances of the adaptive shortest path problem}
     \label{fig:processing_time1}
\end{figure} 

Another important performance measure for comparing algorithms is the processing (CPU) time. We compare the CPU time of our algorithm to those of the iterative algorithm by \cite{chassein2019faster} and the scenario generation approach by \cite{arslan2020min}. Based on preliminary results, the other algorithms were too slow in comparison to the ones selected, thus they were excluded. Figure \ref{fig:processing_time1} illustrate the average CPU times for different algorithms of the instanced solved within the cut-off time for different problem sizes and with \(\Gamma=3\) and \(\Gamma=6\).
It can be seen that the iterative algorithm had better performance than our algorithm in small-size instances (with an average difference of about 2 seconds). However, for the largest problem size ($\mathcal{V}=50$), the iterative algorithm could not solve any instance with $K>3$ within the cut-off time, whereas our approach solved instances of the same size with $K=6$ in less than $20$ seconds. We observe that the double-oracle algorithm has much smaller CPU times for large instances in comparison to the iterative algorithm, which also does not guarantee optimality. The average CPU times for the double-oracle algorithm and the scenario generation algorithm of \cite{arslan2020min} are almost identical. For example, the average CPU times of our algorithm for instances with $\mathcal{V}\in \{20, 25, 40\}$, $\Gamma=3$, and $K=\{2,3,4,5,6\}$  were between $1.55$ and $86.95$ seconds, while those for the scenario generation algorithm were $1.62-82.07$ seconds. However, the scenario generation algorithm could not solve any instances with $\mathcal{V}=50$, and $K>3$, while our algorithm solved some of these instances within the cut-off time.


\subsection{Adaptive Knapsack Problem}

We then tested on the adaptive knapsack problem, for which the objective is to prepare $K$ different combinations of items that respect the capacity constraint without exact knowledge of their profit. Let \(n \in \{1,...,n\}\) be the set of items, \(w_{i}\) and \(p_{i}\), respectively, be the weight and profit of item $i\in n$, and \(b\) is the available budget. The feasible set is defined as \(\mathcal{X}=\{\mathrm{x}\in \{0,1\}^{n}|\;\sum_{i\in N}w_{i}x_{i}\leq b\}\). The goal is to find the best combination of items that maximizes the profit \(\mathrm{p}'\mathrm{x}\). The uncertain parameter \(p_{i}\) is assumed to follow \(p_{i}=(1+\sum_{j\in m}\frac{\Phi_{ij}\xi_{j}}{2})\bar{p_{i}}\), where \(\bar{p}_{i}\) is the nominal profit, \(|m|\) is the number of uncertain factors and \(\Phi \in \mathbb{R}^{|n|\times |m|}\) is the factor loading matrix. The $i$-th row of \(\Phi\) is characterized by the set \(\{\Phi_{i} \in [-1,1]^{|m|}\;|\;\sum_{j\in m}|\Phi_{ij}|=1\}\). As a result, the realized profit of each object \(i \in n\) remains within the interval \([\bar{p}_{i}-\frac{\bar{p}_{i}}{2},\bar{p}_{i}+\frac{\bar{p}_{i}}{2}]\). We solve the problem in different sizes \(n \in \{100,150,200,300\}\) and different values of \(K \in \{2,3,4,5,6\}\). Ten instances of each combination of $n$, and $K$ are solved, and the results obtained from our algorithm are compared to those of the the IA of \cite{chassein2019faster}, the MILP reformulation of \cite{hanasusanto2015k}, and the SG method of \cite{arslan2020min}. These results are summarized in Tables \ref{tab:multicol3} and \ref{tab:multicol4}.

\begin{figure}[H]
    \centering
    \includegraphics[width=\textwidth,height=\textheight,keepaspectratio]{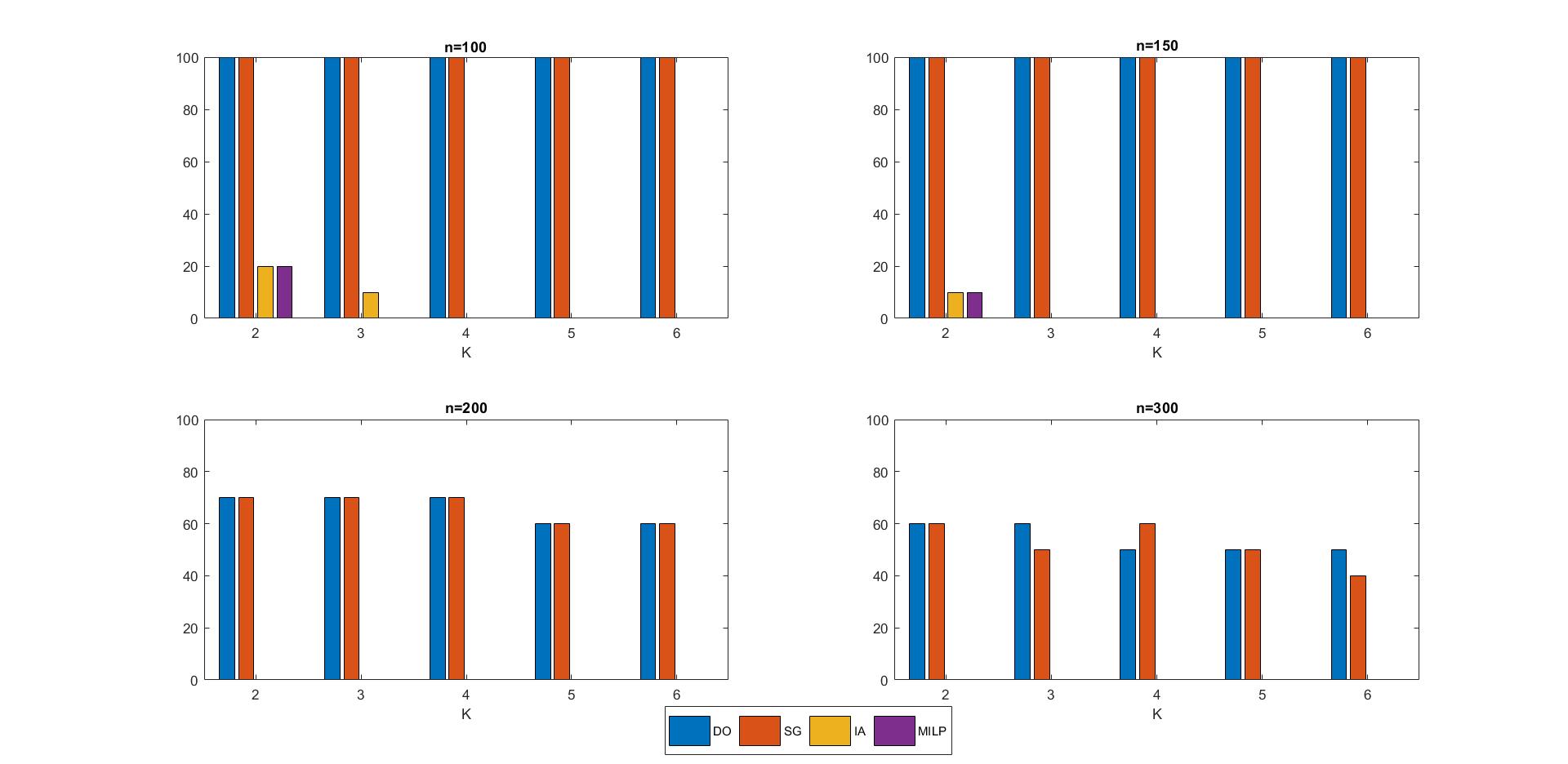}
     \caption{Percentage of solved instances of the adaptive knapsack problem}
     \label{fig:AKP solved problems}
\end{figure}

\begin{figure}[H]
    \centering
    \includegraphics[width=\textwidth,height=\textheight,keepaspectratio]{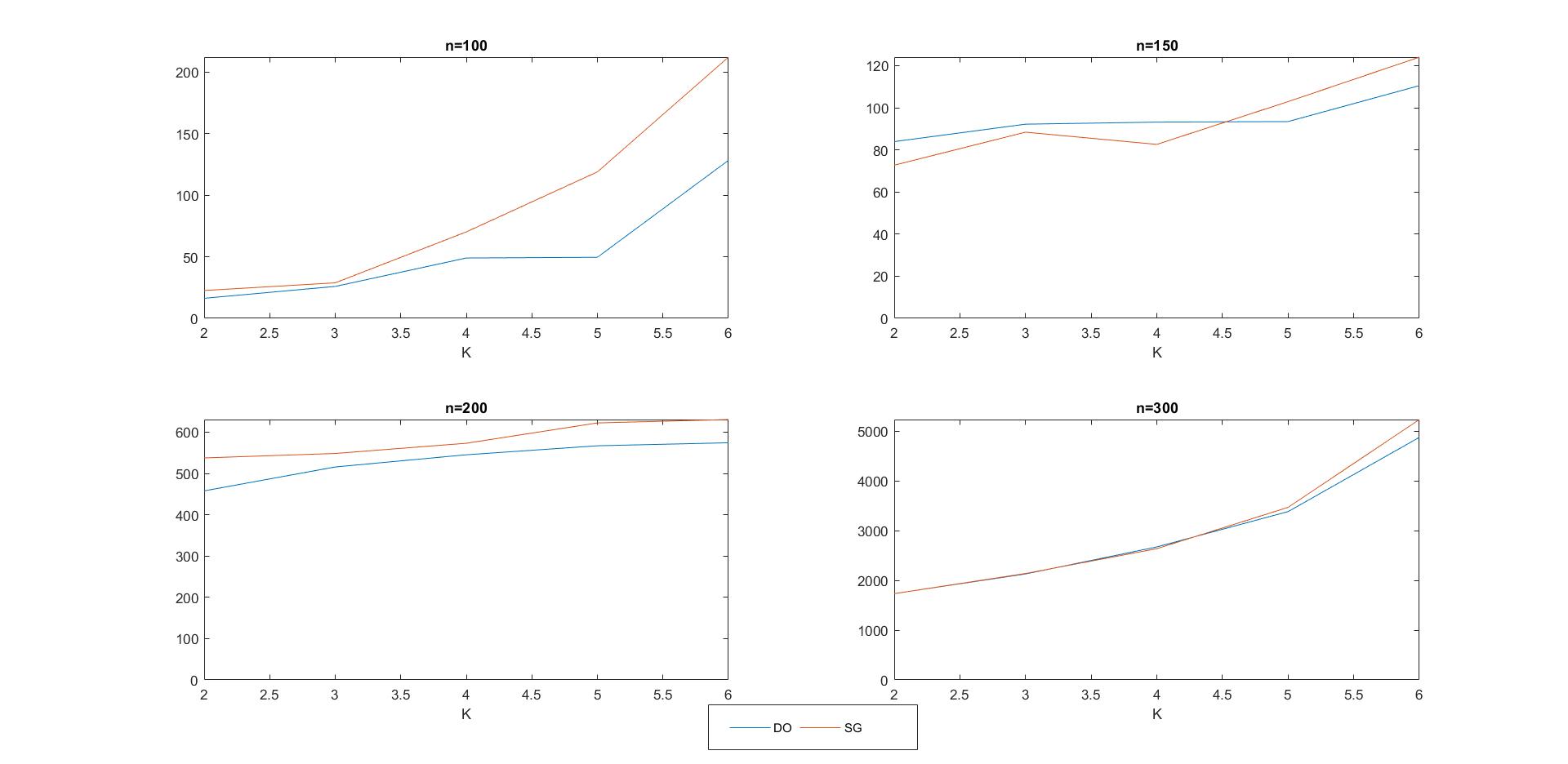}
     \caption{Average CPU times for the solved instances of the adaptive knapsack problem}
     \label{fig:AKP processing time}
\end{figure}

Figure \ref{fig:AKP solved problems} shows the percentage of problems solved to proven optimality by each algorithm. It can be seen that for different problem sizes, our algorithm and the scenario generation method were able to solve almost the same percentage of test instances. However, MILP and IA algorithms could not solve any instances with $n=100$ and $K>2$. By increasing $n$ and $K$, the percentage of solved instances by SG and DO dropped to $50-60\%$, showing the significant impact of $n$ and $K$ on their performance. Figure \ref{fig:AKP processing time} shows the average CPU times for the algorithms with different problem sizes. Since the MILP reformulation and IA could not solve large-size instances, we only compared against the SG algorithm in this round of experiments. It can be seen that the average CPU times of our algorithm and SG approach are very close, with no clear advantage for either algorithm. However, our algorithm is faster than IA and MILP reformulation.
 
\subsection{Problems with Actual First-Stage Decisions}

Finally, we tested on binary $K$-adaptability problems with actual first-stage decisions of the form:    
\begin{subequations}\label{P12}
\begin{alignat}{4}
&\min_{\mathrm{x},\{\mathrm{y}^{k}\}_{k\in [K]}}\max_{\upxi \in \Xi}\min_{k \in [K]}\sum_{i}a_{i}x_{i}+\sum_{j}c_{j}y_{j}^{k}\\
&\;\text{s.t.}\sum_{i}x_{i}=b,\\
&\;\sum_{i}d_{i}x_{i}+\sum_{j}f_{j}y_{j}\geq l,
\end{alignat}
\end{subequations}
where \(\mathrm{x} \in \{0,1\}^n, \mathrm{y} \in \{0,1\}^m\) are the first- and second-stage decision variables, respectively. Moreover, we set \(c_{j}=\bar{c}_{j}-\xi_{j}\hat{c}_{j}\), where \(\bar{c}_{j}\) is the nominal value that is drawn at random from the uniform distribution \(U(8,12)\) and \(\hat{c}_{j}\) is the maximal deviation, set equal to $25\%$ of the nominal value. Moreover, \(a_{i}\), \(d_{i}\) \(f_{j}\) are generated randomly based on the uniform distributions \(U(8,12)\), \(U(50,100)\) and \(U(80,90)\), respectively. Finally, we set \(b=10\), and \(l=0\). The uncertain parameter \(\xi_{j}\) follows \(\Xi=\{\upxi \in [0,1]|\sum_{j}\xi_{j}\leq \Gamma\}\). 10 random instances of each size and uncertainty budget combination were solved. We compared the results of our algorithm to those of the MILP reformulation of \cite{hanasusanto2015k} and Branch-and-Bound (BB) approach of \cite{subramanyam2020k}, which is adapted to each application using the authors’ implementation available \href{https://github.com/AnirudhSubramanyam/KAdaptabilitySolver/blob/v1.0/README.md}{here}, for different instance sizes \(n \in \{20,30,40,50\}\) and \(m \in \{20,30,40,50\}\). Complete results are presented in Tables \ref{tab:multicol5} and \ref{tab:multicol6} in Appendix \ref{AppendixB}.

\begin{figure}[H]
    \centering
    \includegraphics[width=\textwidth,height=\textheight,keepaspectratio]{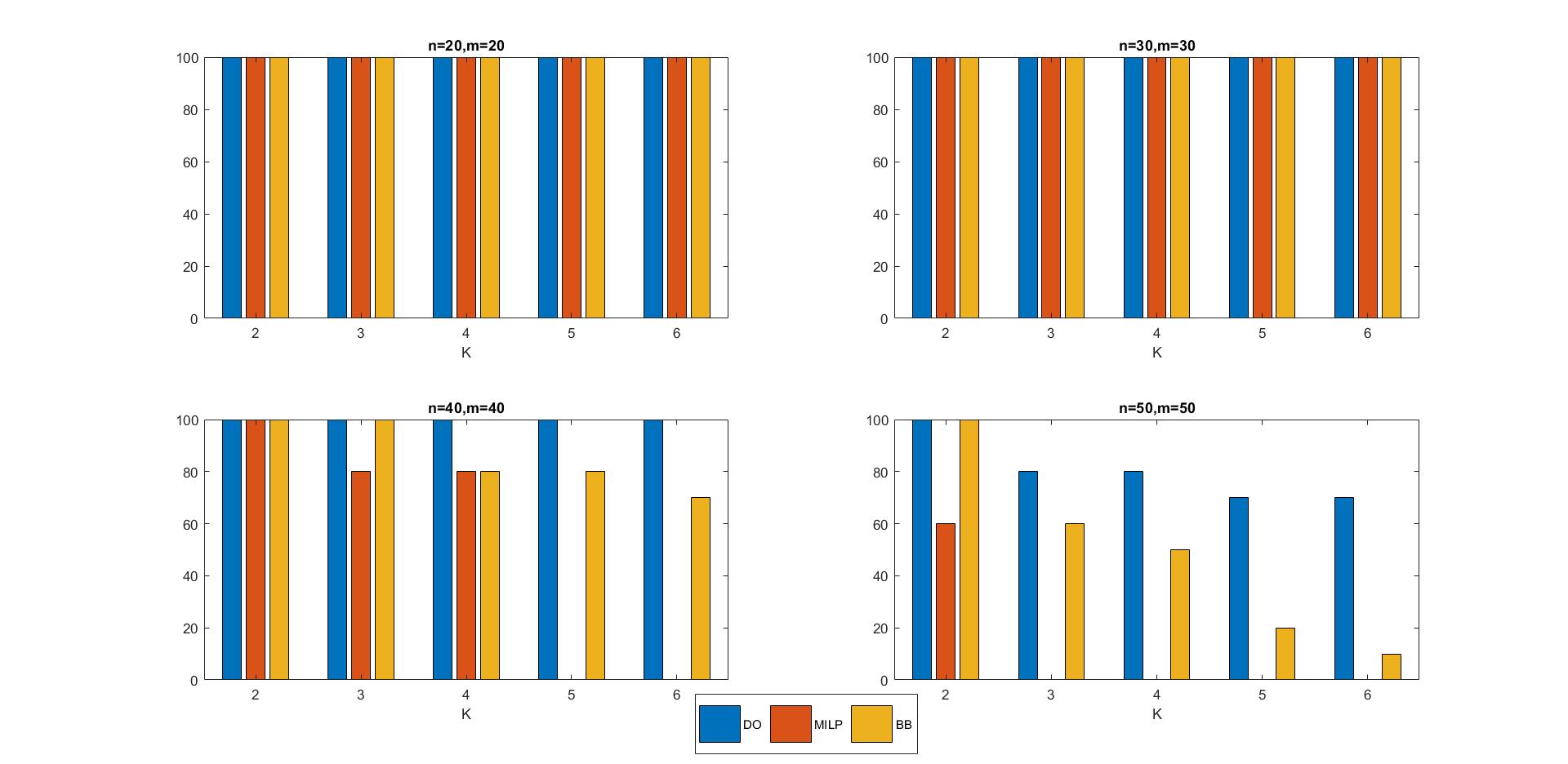}
     \caption{Percentage of solved instances of the two-stage problem}
     \label{fig:two stage percentage}
\end{figure}

\begin{figure}[H]
    \centering
    \includegraphics[width=\textwidth,height=\textheight,keepaspectratio]{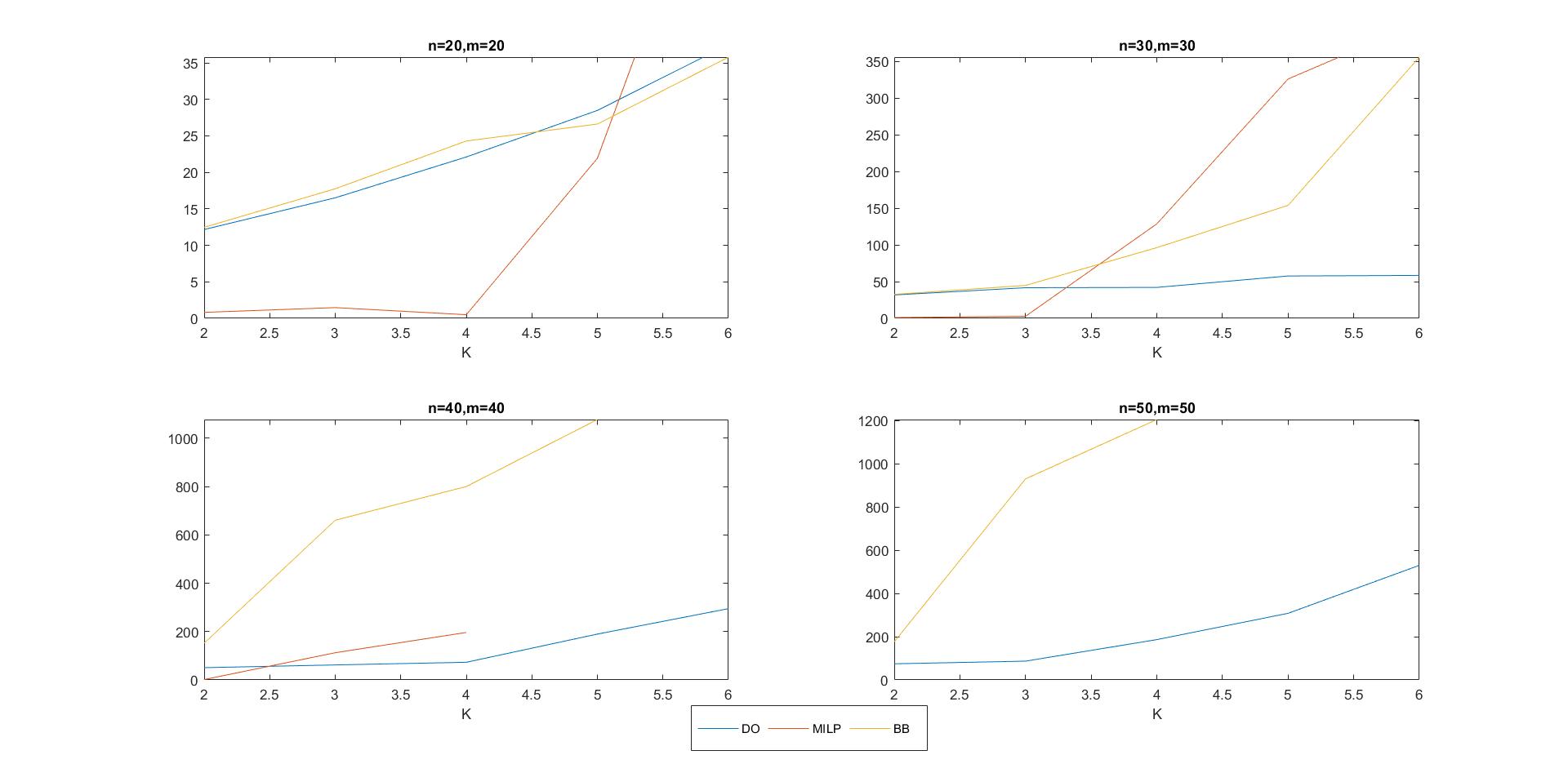}
     \caption{Average CPU times for the solved instances of the two-stage problem}
     \label{fig:Average processing time 2 stage}
\end{figure}

Figure \ref{fig:two stage percentage} shows the percentage of instances solved by each algorithm. For small problem sizes (\ie $n=m=20, 30$), all algorithms were able to solve all instances. However, as the problem sizes were increased, the performance advantage of our algorithm became clear, especially for large values of $K$. For example, the MILP reformulation could solve 60\% of instances of size $n=m=50$ only with $K=2$, but none when $K>2$. In contrast, our algorithm solved the vast majority of instances of the same size to proven optimality within the cut-off time even with $K=6$. Similar insights could be drawn from Figure \ref{fig:Average processing time 2 stage}, which shows the CPU times of all algorithms for the problems that were solved within the cutoff time. Again, it is clear that while the MILP reformulation could solve small instances with small $K$ values efficiently, our algorithm significantly outperforms it in large instances. Furthermore, the proposed algorithm scales well in $K$, thus can result in high adaptability in the face of parameter uncertainty.

\section{Conclusion}\label{sec6}
In this paper, we proposed a new algorithm for solving $K$-adaptability problems, possibly involving constraint and objective function uncertainty and non-linear functions. The algorithm combines discrete scenario generation and $p$-center assignment problem with Logic-based Benders decomposition, where the sub-problem is a MMMRCO. Numerical experiments on benchmark test instances demonstrated that the proposed algorithm performs well with large-size instances and for large $K$ values. It was able to solve larger instances to optimality and produce the result faster than the state-of-the-art algorithms from the existing $K$-adaptability and MMMRCO literature. The proposed algorithm could solve instances of adaptive knapsack problem with up to $300$ items with $k$ up to $6$, while the existing algorithms could solve only a small percentage of the instances with up to $100$ items and $250$ items with $k$ up to $4$. We also showed that the proposed algorithm can solve the shortest path problem with up to $50$ nodes while the existing algorithm barely can solve this problem with more than $50$ nodes. 

The proposed algorithm dominates the iterative scheme of \cite{chassein2019faster} from different points of view. First, our algorithm is faster than the iterative algorithm in large-scale adaptive Knapsack and shortest path problems. Second, The iterative scheme of \cite{chassein2019faster} reaches its limit with $K=2$ with $50$ nodes for the shortest path problem while our algorithm easily solves problems with larger $K$ values up to $6$. Our algorithm also performs better than the RCG algorithm of \cite{goerigk2020min} that is tailored solely for budgeted uncertainty sets. The RCG algorithm still cannot handle the shortest path problem with $n>40$ and $K>3$, whereas results show that our algorithm works efficiently for $n>50$. Our algorithm also dominates the MILP reformulation of \cite{hanasusanto2015k} that cannot solve the shortest path problem with $40$ nodes and $K>3$. Moreover, the MILP reformulation cannot handle the knapsack problem with more than $100$ items. Our algorithm also shows better performance with respect to the percentage of solved instances and processing time than the SG method of \cite{arslan2020min} for the shortest path problem. However, the proposed double-oracle algorithm and the SG method have almost the same performance on the adaptive Knapsack problem. 

Finally, the proposed algorithm leverages a logic-based Benders decomposition to solve the $K$-adaptability faster than any algorithm in the literature. It can handle instances with up to $50$ binary decision variables in each of the first- and second-stage problems. Consequently, we can conclude that not only is our algorithm more efficient with regard to solving large-scale problems in the reasonable time, but also the algorithm is generic and can be extended to tackle different variants of the problems such as nonlinear objective functions and constraints, uncertain parameters in constraints and first-stage objective function.

For future research, we plan to investigate how the proposed algorithm can be extended to general RO problems with combinatorial recourse and to $K$-adaptability problems with continuous first-stage variables. We also propose using efficient methods to solve the $p$-center problem for future research to increase the efficiency of the proposed algorithm.
Another important future research can be improving the lower bound used in optimality cut, a better lower bound can lead to faster convergence to optimal solution in our proposed algorithm.

\bibliographystyle{apa-good}
\bibliography{main}

\begin{thebibliography}{21}
\expandafter\ifx\csname natexlab\endcsname\relax\def\natexlab#1{#1}\fi
\expandafter\ifx\csname url\endcsname\relax
  \def\url#1{{\tt #1}}\fi
\expandafter\ifx\csname urlprefix\endcsname\relax\def\urlprefix{URL }\fi

\bibitem[{Arslan et~al.(2022)Arslan, Poss, \& Silva}]{arslan2020min}
Arslan, A.~N., Poss, M., \& Silva, M. (2022).
\newblock {Min-max-min robust combinatorial optimization with few recourse
  solutions}.
\newblock Working paper or preprint.
\newline\urlprefix\url{https://hal.archives-ouvertes.fr/hal-02939356}

\bibitem[{Ben-Tal et~al.(2004)Ben-Tal, Goryashko, Guslitzer, \&
  Nemirovski}]{ben2004adjustable}
Ben-Tal, A., Goryashko, A., Guslitzer, E., \& Nemirovski, A. (2004).
\newblock Adjustable robust solutions of uncertain linear programs.
\newblock {\em Mathematical programming\/}, {\em 99\/}(2), 351--376.

\bibitem[{Bertsimas et~al.(2011)Bertsimas, Brown, \&
  Caramanis}]{bertsimas2011theory}
Bertsimas, D., Brown, D.~B., \& Caramanis, C. (2011).
\newblock Theory and applications of robust optimization.
\newblock {\em SIAM review\/}, {\em 53\/}(3), 464--501.

\bibitem[{Bertsimas \& Caramanis(2010)}]{bertsimas2010finite}
Bertsimas, D., \& Caramanis, C. (2010).
\newblock Finite adaptability in multistage linear optimization.
\newblock {\em IEEE Transactions on Automatic Control\/}, {\em 55\/}(12),
  2751--2766.

\bibitem[{Bertsimas \& Georghiou(2015)}]{bertsimas2015design}
Bertsimas, D., \& Georghiou, A. (2015).
\newblock Design of near optimal decision rules in multistage adaptive
  mixed-integer optimization.
\newblock {\em Operations Research\/}, {\em 63\/}(3), 610--627.

\bibitem[{Bertsimas \& Georghiou(2018)}]{bertsimas2018binary}
Bertsimas, D., \& Georghiou, A. (2018).
\newblock Binary decision rules for multistage adaptive mixed-integer
  optimization.
\newblock {\em Mathematical Programming\/}, {\em 167\/}(2), 395--433.

\bibitem[{Bertsimas et~al.(2012)Bertsimas, Litvinov, Sun, Zhao, \&
  Zheng}]{bertsimas2012adaptive}
Bertsimas, D., Litvinov, E., Sun, X.~A., Zhao, J., \& Zheng, T. (2012).
\newblock Adaptive robust optimization for the security constrained unit
  commitment problem.
\newblock {\em IEEE transactions on power systems\/}, {\em 28\/}(1), 52--63.

\bibitem[{Buchheim \& Kurtz(2017)}]{buchheim2017min}
Buchheim, C., \& Kurtz, J. (2017).
\newblock Min--max--min robust combinatorial optimization.
\newblock {\em Mathematical Programming\/}, {\em 163\/}(1-2), 1--23.

\bibitem[{Chassein et~al.(2019)Chassein, Goerigk, Kurtz, \&
  Poss}]{chassein2019faster}
Chassein, A., Goerigk, M., Kurtz, J., \& Poss, M. (2019).
\newblock Faster algorithms for min-max-min robustness for combinatorial
  problems with budgeted uncertainty.
\newblock {\em European Journal of Operational Research\/}, {\em 279\/}(2),
  308--319.

\bibitem[{Chen \& Zhang(2009)}]{chen2009uncertain}
Chen, X., \& Zhang, Y. (2009).
\newblock Uncertain linear programs: Extended affinely adjustable robust
  counterparts.
\newblock {\em Operations Research\/}, {\em 57\/}(6), 1469--1482.

\bibitem[{Dhamdhere et~al.(2005)Dhamdhere, Goyal, Ravi, \&
  Singh}]{dhamdhere2005pay}
Dhamdhere, K., Goyal, V., Ravi, R., \& Singh, M. (2005).
\newblock How to pay, come what may: Approximation algorithms for demand-robust
  covering problems.
\newblock In {\em 46th Annual IEEE Symposium on Foundations of Computer Science
  (FOCS'05)\/}, (pp. 367--376). IEEE.

\bibitem[{Georghiou et~al.(2015)Georghiou, Wiesemann, \&
  Kuhn}]{georghiou2015generalized}
Georghiou, A., Wiesemann, W., \& Kuhn, D. (2015).
\newblock Generalized decision rule approximations for stochastic programming
  via liftings.
\newblock {\em Mathematical Programming\/}, {\em 152\/}(1), 301--338.

\bibitem[{Goerigk et~al.(2020)Goerigk, Kurtz, \& Poss}]{goerigk2020min}
Goerigk, M., Kurtz, J., \& Poss, M. (2020).
\newblock Min--max--min robustness for combinatorial problems with discrete
  budgeted uncertainty.
\newblock {\em Discrete Applied Mathematics\/}, {\em 285\/}, 707--725.

\bibitem[{Hanasusanto et~al.(2015)Hanasusanto, Kuhn, \&
  Wiesemann}]{hanasusanto2015k}
Hanasusanto, G.~A., Kuhn, D., \& Wiesemann, W. (2015).
\newblock K-adaptability in two-stage robust binary programming.
\newblock {\em Operations Research\/}, {\em 63\/}(4), 877--891.

\bibitem[{Iancu et~al.(2013)Iancu, Sharma, \&
  Sviridenko}]{iancu2013supermodularity}
Iancu, D.~A., Sharma, M., \& Sviridenko, M. (2013).
\newblock Supermodularity and affine policies in dynamic robust optimization.
\newblock {\em Operations Research\/}, {\em 61\/}(4), 941--956.

\bibitem[{Jiang et~al.(2012)Jiang, Zhang, Li, \& Guan}]{jiang2012benders}
Jiang, R., Zhang, M., Li, G., \& Guan, Y. (2012).
\newblock Benders' decomposition for the two-stage security constrained robust
  unit commitment problem.
\newblock In {\em IIE Annual Conference. Proceedings\/}, (p.~1). Institute of
  Industrial and Systems Engineers (IISE).

\bibitem[{Kuhn et~al.(2011)Kuhn, Wiesemann, \& Georghiou}]{kuhn2011primal}
Kuhn, D., Wiesemann, W., \& Georghiou, A. (2011).
\newblock Primal and dual linear decision rules in stochastic and robust
  optimization.
\newblock {\em Mathematical Programming\/}, {\em 130\/}(1), 177--209.

\bibitem[{Laporte \& Louveaux(1993)}]{laporte1993integer}
Laporte, G., \& Louveaux, F.~V. (1993).
\newblock The integer l-shaped method for stochastic integer programs with
  complete recourse.
\newblock {\em Operations research letters\/}, {\em 13\/}(3), 133--142.

\bibitem[{Subramanyam et~al.(2020)Subramanyam, Gounaris, \&
  Wiesemann}]{subramanyam2020k}
Subramanyam, A., Gounaris, C.~E., \& Wiesemann, W. (2020).
\newblock K-adaptability in two-stage mixed-integer robust optimization.
\newblock {\em Mathematical Programming Computation\/}, {\em 12\/}(2),
  193--224.

\bibitem[{Thiele et~al.(2009)Thiele, Terry, \& Epelman}]{thiele2009robust}
Thiele, A., Terry, T., \& Epelman, M. (2009).
\newblock Robust linear optimization with recourse.
\newblock {\em Rapport technique\/}, (pp. 4--37).

\bibitem[{Zhao \& Zeng(2012)}]{zhao2012robust}
Zhao, L., \& Zeng, B. (2012).
\newblock Robust unit commitment problem with demand response and wind energy.
\newblock In {\em 2012 IEEE power and energy society general meeting\/}, (pp.
  1--8). IEEE.

\end{thebibliography}

\appendix

\section{Proofs}

\begin{proposition}
The objective function value, $w$, of the $p$-center problem \eqref{P5} can be achieved by solving \(w^{*}=\max_{h \in [H]} \min_{j \in [J']}\upxi'_{h}\mathrm{Q}\mathrm{y}_{j}\)
\end{proposition}

\begin{proof}
There are $|H|$ constraints in the form of \eqref{con11} while constraints \eqref{con21} force the problem to select one pair of scenario and solution in each constraint \eqref{con11}. The objective value should be minimum of $w$ that is greater than selected pairs of scenarios and solutions in each $h$ constraint \eqref{con11}. To find the optimal value of $w$, minimum of \(\upxi'_{h}\mathrm{Q}\mathrm{y}_{j}\) for each \(h \in [H']\) is selected, then optimal $w^{*}$ will be maximum of selected pairs of scenarios and solutions in each constraint \eqref{con11}. If \(w^{*}> \max_{h \in [H]} \min_{j \in [J']}\upxi'_{h}\mathrm{Q}\mathrm{y}_{j}\), then it cannot be optimal because there is a feasible solution with lower objective value which is \(\max_{h \in [H]} \min_{j \in [J']}\upxi'_{h}\mathrm{Q}\mathrm{y}_{j}\). On the other hand, if \(w^{*} < \max_{h \in [H]} \min_{j \in [J']}\upxi'_{h}\mathrm{Q}\mathrm{y}_{j}\), then some of constraints \eqref{con11} whose pair of solutions and scenarios are greater than $w$ will be violated. Consequently, \(UB=w^{*}=\max_{h \in [H]} \min_{j \in [J']}\upxi'_{h}\mathrm{Q}\mathrm{y}_{j}\).  
\end{proof}

\begin{proposition}
The objective function value, $\gamma$, of the problem \eqref{P7} can be achieved by solving \(\gamma^{*}=\max_{h \in [H']} \min_{\mathrm{y}^{k}\in \mathcal{Y}}\upxi'_h\mathrm{Q}\mathrm{y}^{k}\).
\end{proposition}

\begin{proof}
Let us assume that \(\mathrm{y}^{{k}^{*}}\) is the optimal solution of problem \eqref{P7}. There are $|K|\times |H'|$ constraints of \eqref{con62}. For each $h\in [H']$ there are $|K|$ constraints in form of \eqref{con62}. Constraints \eqref{con63} force the problem to select one pair of scenarios and solutions for each \(h \in [H']\). The objective function is minimization, consequently, minimum of $(\mathrm{y}^{{k}^{*}},\upxi'_{h})$ for each $h\in [H']$ are selected. Since $\gamma$ should be greater than \(\upxi'_{h}\mathrm{Q}\mathrm{y}^{{k}^{*}}\), then maximum of selected pairs will be optimal objective value $\gamma^*$. Consequently, the optimal value of objective function, $\gamma^{*}$ can be achieved by solving \(\gamma^{*}=\max_{h \in [H']} \min_{\mathrm{y}^{k}\in \mathcal{Y}}\upxi'_h\mathrm{Q}\mathrm{y}^{k}\).
\end{proof}

\section{Results}\label{AppendixB}

\begin{table}[ht]
\caption{Percentage of solved shortest path instances for $\Gamma=3$}
\begin{center}
\begin{tabular}{|c|l|c|c|c|c|c|}
    \hline
    Size&Algorithms&$k=2$&$k=3$&$k=4$&$k=5$&$k=6$\\
    \hline
    \hline
    \multirow{5}{*}{20}&Double-Oracle &$100\%$&$100\%$&$100\%$&$100\%$&$100\%$\\
    &MILP \cite{hanasusanto2015k}&$100\%$&$100\%$&$80\%$&$80\%$&$80\%$\\
    &IA \cite{chassein2019faster}&$90\%$&$90\%$&$70\%$&$70\%$&$60\%$\\
    &RCG \cite{goerigk2020min}&$40\%$&$30\%$&$10\%$&$10\%$&$10\%$\\
    &SG \cite{arslan2020min}&$100\%$&$100\%$&$80\%$&$80\%$&$80\%$\\
    \hline
    \multirow{5}{*}{25}&Double-Oracle &$100\%$&$100\%$&$100\%$&$100\%$&$100\%$\\
    &MILP \cite{hanasusanto2015k}&$100\%$&$100\%$&$40\%$&$10\%$&$10\%$\\
    &IA \cite{chassein2019faster}&$100\%$&$100\%$&$80\%$&$40\%$&$40\%$\\
    &RCG \cite{goerigk2020min}&$50\%$&$20\%$&$10\%$&$0\%$&$0\%$\\
    &SG \cite{arslan2020min}&$100\%$&$100\%$&$100\%$&$100\%$&$100\%$\\
    \hline
    \multirow{5}{*}{40}&Double-Oracle &$100\%$&$100\%$&$100\%$&$100\%$&$100\%$\\
    &MILP \cite{hanasusanto2015k}&$90\%$&$80\%$&$40\%$&$0\%$&$0\%$\\
    &IA \cite{chassein2019faster}&$100\%$&$100\%$&$100\%$&$80\%$&$70\%$\\
    &RCG \cite{goerigk2020min}&$10\%$&$10\%$&$0\%$&$0\%$&$0\%$\\
    &SG \cite{arslan2020min}&$100\%$&$100\%$&$80\%$&$80\%$&$80\%$\\
    \hline
    \multirow{5}{*}{50}&Double-Oracle &$60\%$&$40\%$&$30\%$&$30\%$&$30\%$\\
    &MILP \cite{hanasusanto2015k}&$10\%$&$0\%$&$0\%$&$0\%$&$0\%$\\
    &IA \cite{chassein2019faster}&$20\%$&$10\%$&$0\%$&$0\%$&$0\%$\\
    &RCG \cite{goerigk2020min}&$0\%$&$0\%$&$0\%$&$0\%$&$0\%$\\
    &SG \cite{arslan2020min}&$30\%$&$10\%$&$0\%$&$0\%$&$0\%$\\
    \hline
\end{tabular}
\end{center}
\label{tab:multicol1}
\end{table}

\begin{table}[ht]
\caption{Percentage of solved shortest path instances for $\Gamma=6$}
\begin{center}
\begin{tabular}{|c|l|c|c|c|c|c|}
    \hline
    Size&Algorithms&$k=2$&$k=3$&$k=4$&$k=5$&$k=6$\\
    \hline
    \hline
\multirow{5}{*}{20}	&	Double-Oracle	&	$	100	\%	$	&	$	100	\%	$	&	$	100	\%	$	&	$	100	\%	$	&	$	90	\%	$	\\
	&	MILP \cite{hanasusanto2015k}	&	$	100	\%	$	&	$	100	\%	$	&	$	80	\%	$	&	$	70	\%	$	&	$	70	\%	$	\\
	&	IA \cite{chassein2019faster}	&	$	90	\%	$	&	$	80	\%	$	&	$	70	\%	$	&	$	60	\%	$	&	$	60	\%	$	\\
	&	RCG \cite{goerigk2020min}	&	$	40	\%	$	&	$	30	\%	$	&	$	10	\%	$	&	$	0	\%	$	&	$	0	\%	$	\\
	&	SG \cite{arslan2020min}	&	$	80	\%	$	&	$	80	\%	$	&	$	60	\%	$	&	$	60	\%	$	&	$	50	\%	$	\\
\hline																												
\multirow{5}{*}{25}	&	Double-Oracle	&	$	100	\%	$	&	$	100	\%	$	&	$	100	\%	$	&	$	100	\%	$	&	$	100	\%	$	\\
	&	MILP \cite{hanasusanto2015k}	&	$	100	\%	$	&	$	100	\%	$	&	$	40	\%	$	&	$	10	\%	$	&	$	0	\%	$	\\
	&	IA \cite{chassein2019faster}	&	$	100	\%	$	&	$	100	\%	$	&	$	100	\%	$	&	$	60	\%	$	&	$	40	\%	$	\\
	&	RCG \cite{goerigk2020min}	&	$	40	\%	$	&	$	10	\%	$	&	$	0	\%	$	&	$	0	\%	$	&	$	0	\%	$	\\
	&	SG \cite{arslan2020min}	&	$	90	\%	$	&	$	90	\%	$	&	$	70	\%	$	&	$	70	\%	$	&	$	70	\%	$	\\
\hline																												
\multirow{5}{*}{40}	&	Double-Oracle	&	$	90	\%	$	&	$	90	\%	$	&	$	80	\%	$	&	$	80	\%	$	&	$	80	\%	$	\\
	&	MILP \cite{hanasusanto2015k}	&	$	80	\%	$	&	$	70	\%	$	&	$	10	\%	$	&	$	0	\%	$	&	$	0	\%	$	\\
	&	IA \cite{chassein2019faster}	&	$	100	\%	$	&	$	100	\%	$	&	$	100	\%	$	&	$	90	\%	$	&	$	50	\%	$	\\
	&	RCG \cite{goerigk2020min}	&	$	10	\%	$	&	$	0	\%	$	&	$	0	\%	$	&	$	0	\%	$	&	$	0	\%	$	\\
	&	SG \cite{arslan2020min}	&	$	80	\%	$	&	$	80	\%	$	&	$	60	\%	$	&	$	50	\%	$	&	$	50	\%	$	\\
\hline																												
\multirow{5}{*}{50}	&	Double-Oracle	&	$	50	\%	$	&	$	40	\%	$	&	$	10	\%	$	&	$	10	\%	$	&	$	10	\%	$	\\
	&	MILP \cite{hanasusanto2015k}	&	$	0	\%	$	&	$	0	\%	$	&	$	0	\%	$	&	$	0	\%	$	&	$	0	\%	$	\\
	&	IA \cite{chassein2019faster}	&	$	10	\%	$	&	$	0	\%	$	&	$	0	\%	$	&	$	0	\%	$	&	$	0	\%	$	\\
	&	RCG \cite{goerigk2020min}	&	$	0	\%	$	&	$	0	\%	$	&	$	0	\%	$	&	$	0	\%	$	&	$	0	\%	$	\\
	&	SG \cite{arslan2020min}	&	$	30	\%	$	&	$	10	\%	$	&	$	0	\%	$	&	$	0	\%	$	&	$	0	\%	$	\\
    \hline
\end{tabular}
\end{center}
\label{tab:multicol2}
\end{table}

\begin{table}[ht]
\caption{Average CPU time for solved shortest path instances for $\Gamma=3$}
\begin{center}
\begin{tabular}{|c|l|c|c|c|c|c|}
    \hline
    Size&Algorithms&$k=2$&$k=3$&$k=4$&$k=5$&$k=6$\\
    \hline
    \hline
\multirow{5}{*}{20}	&	Double-Oracle	&	$	1.55	$	&	$	1.56	$	&	$	1.69	$	&	$	1.82	$	&	$	1.99	$	\\
	&	MILP \cite{hanasusanto2015k}	&	$	0.46	$	&	$	0.96	$	&	$	3.09	$	&	$	6.54	$	&	$	433.08	$	\\
	&	IA \cite{chassein2019faster}	&	$	0.35	$	&	$	0.41	$	&	$	0.38	$	&	$	0.40	$	&	$	0.43	$	\\
	&	RCG \cite{goerigk2020min}	&	$	1749.49	$	&	$	2569.55	$	&	$	2955.92	$	&	$	3491.93	$	&	$	3564.87	$	\\
	&	SG \cite{arslan2020min}	&	$	1.62	$	&	$	1.59	$	&	$	1.77	$	&	$	1.85	$	&	$	1.99	$	\\
\hline																							
\multirow{5}{*}{25}	&	Double-Oracle	&	$	2.87	$	&	$	2.80	$	&	$	3.26	$	&	$	3.70	$	&	$	7.23	$	\\
	&	MILP \cite{hanasusanto2015k}	&	$	12.69	$	&	$	53.75	$	&	$	524.45	$	&	$	2.50	$	&	$	16.55	$	\\
	&	IA \cite{chassein2019faster}	&	$	3.44	$	&	$	3.42	$	&	$	3.23	$	&	$	4.08	$	&	$	4.36	$	\\
	&	RCG \cite{goerigk2020min}	&	$	4373.73	$	&	$	6423.87	$	&	$	7129.75	$	&		NA		&		NA		\\
	&	SG \cite{arslan2020min}	&	$	2.93	$	&	$	2.95	$	&	$	3.07	$	&	$	3.94	$	&	$	6.52	$	\\
\hline																							
\multirow{5}{*}{40}	&	Double-Oracle	&	$	8.63	$	&	$	7.50	$	&	$	9.04	$	&	$	11.97	$	&	$	86.95	$	\\
	&	MILP \cite{hanasusanto2015k}	&	$	90.61	$	&	$	299.26	$	&	$	873.76	$	&		NA		&		NA		\\
	&	IA \cite{chassein2019faster}	&	$	68.45	$	&	$	81.99	$	&	$	79.13	$	&	$	65.68	$	&	$	77.08	$	\\
	&	RCG \cite{goerigk2020min}	&	$	3498.99	$	&	$	5139.10	$	&		NA		&		NA		&		NA		\\
	&	SG \cite{arslan2020min}	&	$	6.91	$	&	$	6.51	$	&	$	8.93	$	&	$	10.85	$	&	$	82.07	$	\\
\hline																							
\multirow{5}{*}{50}	&	Double-Oracle	&	$	9.00	$	&	$	17.04	$	&	$	17.15	$	&	$	15.36	$	&	$	17.78	$	\\
	&	MILP \cite{hanasusanto2015k}	&	$	355.99	$	&		NA		&		NA		&		NA		&		NA		\\
	&	IA \cite{chassein2019faster}	&	$	102.68	$	&	$	118.69	$	&		NA		&		NA		&		NA		\\
	&	RCG \cite{goerigk2020min}	&		NA		&		NA		&		NA		&		NA		&		NA		\\
	&	SG \cite{arslan2020min}	&	$	28.93	$	&		NA		&		NA		&		NA		&		NA		\\
    \hline
\end{tabular}
\end{center}
\label{tab:multicol3}
\end{table}

\begin{table}[ht]
\caption{Average CPU time (s) for solved shortest path instances for $\Gamma=6$}
\begin{center}
\begin{tabular}{|c|l|c|c|c|c|c|}
    \hline
    Size&Algorithms&$k=2$&$k=3$&$k=4$&$k=5$&$k=6$\\
    \hline
    \hline
\multirow{5}{*}{20}	&	Double-Oracle	&	$	1.54	$	&	$	1.54	$	&	$	1.54	$	&	$	1.61	$	&	$	1.85	$	\\
	&	MILP \cite{hanasusanto2015k}	&	$	0.72	$	&	$	1.26	$	&	$	3.31	$	&	$	7.21	$	&	$	119.40	$	\\
	&	IA \cite{chassein2019faster}	&	$	0.28	$	&	$	0.35	$	&	$	0.36	$	&	$	0.40	$	&	$	0.42	$	\\
	&	RCG \cite{goerigk2020min}	&	$	1766.99	$	&	$	2595.24	$	&	$	2985.48	$	&		NA		&		NA		\\
	&	SG \cite{arslan2020min}	&	$	1.55	$	&	$	1.57	$	&	$	1.61	$	&	$	1.64	$	&	$	1.93	$	\\
\hline																							
\multirow{5}{*}{25}	&	Double-Oracle	&	$	4.43	$	&	$	4.47	$	&	$	5.08	$	&	$	8.19	$	&	$	9.23	$	\\
	&	MILP \cite{hanasusanto2015k}	&	$	14.31	$	&	$	75.40	$	&	$	151.24	$	&	$	389.98	$	&		NA		\\
	&	IA \cite{chassein2019faster}	&	$	3.39	$	&	$	4.34	$	&	$	4.84	$	&	$	5.28	$	&	$	14.11	$	\\
	&	RCG \cite{goerigk2020min}	&	$	4417.47	$	&	$	6488.11	$	&		NA		&		NA		&		NA		\\
	&	SG \cite{arslan2020min}	&	$	4.40	$	&	$	4.66	$	&	$	4.91	$	&	$	8.16	$	&	$	10.10	$	\\
\hline																							
\multirow{5}{*}{40}	&	Double-Oracle	&	$	10.59	$	&	$	11.27	$	&	$	14.58	$	&	$	22.20	$	&	$	26.96	$	\\
	&	MILP \cite{hanasusanto2015k}	&	$	70.59	$	&	$	664.92	$	&	$	1044.60	$	&		NA		&		NA		\\
	&	IA \cite{chassein2019faster}	&	$	54.03	$	&	$	73.24	$	&	$	74.91	$	&	$	81.40	$	&	$	85.20	$	\\
	&	RCG \cite{goerigk2020min}	&	$	3533.98	$	&		NA		&		NA		&		NA		&		NA		\\
	&	SG \cite{arslan2020min}	&	$	11.41	$	&	$	11.81	$	&	$	12.20	$	&	$	20.76	$	&	$	24.96	$	\\
\hline																							
\multirow{5}{*}{50}	&	Double-Oracle	&	$	10.46	$	&	$	15.45	$	&	$	15.70	$	&	$	18.01	$	&	$	22.08	$	\\
	&	MILP \cite{hanasusanto2015k}	&		NA		&		NA		&		NA		&		NA		&		NA		\\
	&	IA \cite{chassein2019faster}	&	$	81.04	$	&		NA		&		NA		&		NA		&		NA		\\
	&	RCG \cite{goerigk2020min}	&		NA		&		NA		&		NA		&		NA		&		NA		\\
	&	SG \cite{arslan2020min}	&	$	33.67	$	&	$	48.21	$	&		NA		&		NA		&		NA		\\
    \hline
\end{tabular}
\end{center}
\label{tab:multicol4}
\end{table}

\begin{table}[ht]
\caption{Percentage of solved knapsack problem instances}
\begin{center}
\begin{tabular}{|c|l|c|c|c|c|c|}
    \hline
    Size&Algorithms&$k=2$&$k=3$&$k=4$&$k=5$&$k=6$\\
    \hline
    \hline
\multirow{4}{*}{100}	&	Double-Oracle	&	$	100	\%	$	&	$	100	\%	$	&	$	100	\%	$	&	$	100	\%	$	&	$	100	\%	$	\\
	&	MILP \cite{hanasusanto2015k}	&	$	20	\%	$	&	$	0	\%	$	&	$	0	\%	$	&	$	0	\%	$	&	$	0	\%	$	\\
	&	IA \cite{chassein2019faster}	&	$	20	\%	$	&	$	10	\%	$	&	$	0	\%	$	&	$	0	\%	$	&	$	0	\%	$	\\
	&	SG \cite{arslan2020min}	&	$	100	\%	$	&	$	100	\%	$	&	$	100	\%	$	&	$	100	\%	$	&	$	100	\%	$	\\
\hline				
\multirow{4}{*}{150}	&	Double-Oracle	&	$	100	\%	$	&	$	100	\%	$	&	$	100	\%	$	&	$	100	\%	$	&	$	100	\%	$	\\
	&	MILP \cite{hanasusanto2015k}	&	$	10	\%	$	&	$	0	\%	$	&	$	0	\%	$	&	$	0	\%	$	&	$	0	\%	$	\\
	&	IA \cite{chassein2019faster}	&	$	10	\%	$	&	$	0	\%	$	&	$	0	\%	$	&	$	0	\%	$	&	$	0	\%	$	\\
	&	SG \cite{arslan2020min}	&	$	100	\%	$	&	$	100	\%	$	&	$	100	\%	$	&	$	100	\%	$	&	$	100	\%	$	\\
\hline																								
\multirow{4}{*}{200}	&	Double-Oracle	&	$	70	\%	$	&	$	70	\%	$	&	$	70	\%	$	&	$	60	\%	$	&	$	60	\%	$	\\
	&	MILP \cite{hanasusanto2015k}	&	$	0	\%	$	&	$	0	\%	$	&	$	0	\%	$	&	$	0	\%	$	&	$	0	\%	$	\\
	&	IA \cite{chassein2019faster}	&	$	0	\%	$	&	$	0	\%	$	&	$	0	\%	$	&	$	0	\%	$	&	$	0	\%	$	\\
	&	SG \cite{arslan2020min}	&	$	70	\%	$	&	$	70	\%	$	&	$	70	\%	$	&	$	60	\%	$	&	$	60	\%	$	\\
\hline																								
\multirow{4}{*}{300}	&	Double-Oracle	&	$	60	\%	$	&	$	60	\%	$	&	$	50	\%	$	&	$	50	\%	$	&	$	50	\%	$	\\
	&	MILP \cite{hanasusanto2015k}	&	$	0	\%	$	&	$	0	\%	$	&	$	0	\%	$	&	$	0	\%	$	&	$	0	\%	$	\\
	&	IA \cite{chassein2019faster}	&	$	0	\%	$	&	$	0	\%	$	&	$	0	\%	$	&	$	0	\%	$	&	$	0	\%	$	\\
	&	SG \cite{arslan2020min}	&	$	60	\%	$	&	$	50	\%	$	&	$	60	\%	$	&	$	50	\%	$	&	$	40	\%	$	\\
    \hline
\end{tabular}
\end{center}
\label{tab:multicol5}
\end{table}

\begin{table}[ht]
\caption{Average CPU time (s) of solved knapsack problem instances}
\begin{center}
\begin{tabular}{|c|l|c|c|c|c|c|}
    \hline
    Size&Algorithms&$k=2$&$k=3$&$k=4$&$k=5$&$k=6$\\
    \hline
    \hline
\multirow{5}{*}{100}	&	Double-Oracle	&	$	16.29	$	&	$	25.92	$	&	$	49.01	$	&	$	49.60	$	&	$	128.35	$	\\
	&	MILP \cite{hanasusanto2015k}	&	$	3829.95	$	&		NA		&		NA		&		NA		&		NA		\\
	&	IA \cite{chassein2019faster}	&	$	999.12	$	&	$	1460.23	$	&		NA		&		NA		&		NA		\\
	&	SG \cite{arslan2020min}	&	$	22.60	$	&	$	28.85	$	&	$	70.23	$	&	$	118.96	$	&	$	212.19	$	\\
\hline																							
\multirow{5}{*}{150}	&	Double-Oracle	&	$	83.93	$	&	$	92.21	$	&	$	93.19	$	&	$	93.45	$	&	$	110.52	$	\\
	&	MILP \cite{hanasusanto2015k}	&	$	4072.81	$	&		NA		&		NA		&		NA		&		NA		\\
	&	IA \cite{chassein2019faster}	&	$	70.94	$	&		NA		&		NA		&		NA		&		NA		\\
	&	SG \cite{arslan2020min}	&	$	72.73	$	&	$	88.41	$	&	$	82.61	$	&	$	102.92	$	&	$	124.05	$	\\
\hline																							
\multirow{5}{*}{200}	&	Double-Oracle	&	$	457.84	$	&	$	515.72	$	&	$	545.24	$	&	$	567.05	$	&	$	574.33	$	\\
	&	MILP \cite{hanasusanto2015k}	&		NA		&		NA		&		NA		&		NA		&		NA		\\
	&	IA \cite{chassein2019faster}	&		NA		&		NA		&		NA		&		NA		&		NA		\\
	&	SG \cite{arslan2020min}	&	$	537.63	$	&	$	548.42	$	&	$	573.13	$	&	$	622.36	$	&	$	630.55	$	\\
\hline																							
\multirow{5}{*}{300}	&	Double-Oracle	&	$	1734.74	$	&	$	2129.96	$	&	$	2673.10	$	&	$	3382.99	$	&	$	4875.86	$	\\
	&	MILP \cite{hanasusanto2015k}	&		NA		&		NA		&		NA		&		NA		&		NA		\\
	&	IA \cite{chassein2019faster}	&		NA		&		NA		&		NA		&		NA		&		NA		\\
	&	SG \cite{arslan2020min}	&	$	1733.66	$	&	$	2138.93	$	&	$	2638.38	$	&	$	3468.43	$	&	$	5233.18	$	\\

    \hline
\end{tabular}
\end{center}
\label{tab:multicol6}
\end{table}

\begin{table}[ht]
\caption{Percentage of solved generic $K$-adaptability instances}
\begin{center}
\begin{tabular}{|c|l|c|c|c|c|c|}
    \hline
    Size&Algorithms&$k=2$&$k=3$&$k=4$&$k=5$&$k=6$\\
    \hline
    \hline
\multirow{3}{*}{n=20,m=20}	&	Double-Oracle	&	$	100	\%	$	&	$	100	\%	$	&	$	100	\%	$	&	$	100	\%	$	&	$	100	\%	$	\\
	&	MILP \cite{hanasusanto2015k}	&	$	100	\%	$	&	$	100	\%	$	&	$	100	\%	$	&	$	100	\%	$	&	$	100	\%	$	\\
	&	BB \cite{subramanyam2020k}	&	$	100	\%	$	&	$	100	\%	$	&	$	100	\%	$	&	$	100	\%	$	&	$	100	\%	$	\\
\hline																												
\multirow{3}{*}{n=30,m=30}	&	Double-Oracle	&	$	100	\%	$	&	$	100	\%	$	&	$	100	\%	$	&	$	100	\%	$	&	$	100	\%	$	\\
	&	MILP \cite{hanasusanto2015k}	&	$	100	\%	$	&	$	100	\%	$	&	$	100	\%	$	&	$	100	\%	$	&	$	100	\%	$	\\
	&	BB \cite{subramanyam2020k}	&	$	100	\%	$	&	$	100	\%	$	&	$	100	\%	$	&	$	100	\%	$	&	$	100	\%	$	\\
\hline																												
\multirow{3}{*}{n=40,m=40}	&	Double-Oracle	&	$	100	\%	$	&	$	100	\%	$	&	$	100	\%	$	&	$	100	\%	$	&	$	100	\%	$	\\
	&	MILP \cite{hanasusanto2015k}	&	$	100	\%	$	&	$	80	\%	$	&	$	80	\%	$	&	$	0	\%	$	&	$	0	\%	$	\\
	&	BB \cite{subramanyam2020k}	&	$	100	\%	$	&	$	100	\%	$	&	$	80	\%	$	&	$	80	\%	$	&	$	70	\%	$	\\
\hline																												
\multirow{3}{*}{n=50,m=50}	&	Double-Oracle	&	$	100	\%	$	&	$	80	\%	$	&	$	80	\%	$	&	$	70	\%	$	&	$	70	\%	$	\\
	&	MILP \cite{hanasusanto2015k}	&	$	60	\%	$	&	$	0	\%	$	&	$	0	\%	$	&	$	0	\%	$	&	$	0	\%	$	\\
	&	BB \cite{subramanyam2020k}	&	$	100	\%	$	&	$	60	\%	$	&	$	50	\%	$	&	$	20	\%	$	&	$	10	\%	$	\\
    \hline
\end{tabular}
\end{center}
\label{tab:multicol7}
\end{table}

\begin{table}[ht]
\caption{Average CPU time (s) of solved generic $K$-adaptability instances}
\begin{center}
\begin{tabular}{|c|l|c|c|c|c|c|}
    \hline
    Size&Algorithms&$k=2$&$k=3$&$k=4$&$k=5$&$k=6$\\
    \hline
    \hline
\multirow{3}{*}{n=20,m=20}	&	Double-Oracle	&	$	12.17	$	&	$	16.51	$	&	$	22.11	$	&	$	28.48	$	&	$	37.52	$	\\
	&	MILP \cite{hanasusanto2015k}	&	$	0.82	$	&	$	1.47	$	&	$	0.50	$	&	$	21.93	$	&	$	70.54	$	\\
	&	BB \cite{subramanyam2020k}	&	$	12.47	$	&	$	17.75	$	&	$	24.29	$	&	$	26.62	$	&	$	35.78	$	\\
\hline																							
\multirow{3}{*}{n=30,m=30}	&	Double-Oracle	&	$	31.84	$	&	$	41.69	$	&	$	42.13	$	&	$	57.68	$	&	$	58.45	$	\\
	&	MILP \cite{hanasusanto2015k}	&	$	1.10	$	&	$	2.78	$	&	$	128.62	$	&	$	325.99	$	&	$	402.52	$	\\
	&	BB \cite{subramanyam2020k}	&	$	32.62	$	&	$	44.82	$	&	$	96.30	$	&	$	153.90	$	&	$	355.75	$	\\
\hline																							
\multirow{3}{*}{n=40,m=40}	&	Double-Oracle	&	$	50.67	$	&	$	62.05	$	&	$	73.02	$	&	$	189.71	$	&	$	294.91	$	\\
	&	MILP \cite{hanasusanto2015k}	&	$	1.50	$	&	$	112.55	$	&	$	196.46	$	&		NA		&		NA		\\
	&	BB \cite{subramanyam2020k}	&	$	81.91	$	&	$	660.71  $	&	$	800.25	$	&	$	1077.28	$	&	NA	\\
\hline																							
\multirow{3}{*}{n=50,m=50}	&	Double-Oracle	&	$	74.63	$	&	$	86.93	$	&	$	186.54	$	&	$	308.29	$	&	$	530.35	$	\\
	&	MILP \cite{hanasusanto2015k}	&	$	3.87	$	&		NA		&		NA		&		NA		&		NA		\\
	&	BB \cite{subramanyam2020k}	&	$	176.47	$	&	$	930.46  $	&	$	1205.01	$	&	NA	&	NA	\\
\hline
\end{tabular}
\end{center}
\label{tab:multicol8}
\end{table}

\end{document}